\documentclass[fleqn,10pt]{wlscirep}
\title{Learning to predict synchronization of coupled oscillators on randomly generated graphs}
\usepackage{amssymb,mathrsfs}
\usepackage{amsmath}
\usepackage{adjustbox}
\usepackage{calc}
\usepackage{array}
\usepackage{tikz}
\usepackage{algorithm}
\usepackage[noend]{algpseudocode}
\usepackage[normalem]{ulem}
\usepackage{array}
\usepackage{enumerate}
\usepackage[english]{babel}
\usepackage{multirow}
\usepackage{float}
\usepackage[caption = false]{subfig}
\usepackage{enumitem}
\usepackage{accents}
\usepackage{subfig}

\usepackage{theoremref}
\newcommand\dataset[1]{\textsc{\texttt{#1}}}
\makeatletter
\renewcommand{\@seccntformat}[1]{}

\makeatother

\numberwithin{theorem}{section}

\usepackage{pifont}
\usepackage{wrapfig}

\def\E{\mathbb{E}}

\makeatletter
\def\namedlabel#1#2{\begingroup
   \def\@currentlabel{#2}%
   \label{#1}\endgroup
}
\makeatother

\newtheorem{innercustomquestion}{}
\newenvironment{customquestion}[1]
{\renewcommand\theinnercustomquestion{#1}\innercustomquestion}
{\endinnercustomquestion}

\definecolor{hancolor}{rgb}{0.0 0.0, 1.0}
\newcommand{\commHL}[1]{{\textcolor{black}{#1}}} 

\author[1,$\dagger$]{Hardeep Bassi} 
\author[2,$\dagger$]{Richard P. Yim} 
\author[3]{Joshua Vendrow} 
\author[3]{Rohith Koduluka}
\author[4]{Cherlin Zhu} 
\author[5,*]{Hanbaek Lyu}
\affil[1]{Department of Applied Mathematics, University of California, Merced, CA 95343, USA}
\affil[2]{Department of Mathematics, University of California, Davis, CA 95616, USA}
\affil[3]{Department of Mathematics, University of California, Los Angeles, CA 90095, USA} 
\affil[4]{
Department of Applied Mathematics and Statistics, Johns Hopkins University, Baltimore, MD 21218, USA}
\affil[5]{
Department of Mathematics, University of Wisconsin - Madison, WI 53706, USA}
\affil[$\dagger$]{Co-first authors}
\affil[*]{Corresponding author; \texttt{hlyu@math.wisc.edu}}
\begin{abstract}
Suppose we are given a system of coupled oscillators on an unknown graph along with the trajectory of the system during some period. Can we predict whether the system will eventually synchronize? Even with a known underlying graph structure, this is an important yet analytically intractable question in general. In this work, we take an \commHL{alternative approach to the synchronization prediction problem by viewing it as a classification problem based on the fact that any given system will eventually synchronize or converge to a non-synchronizing limit cycle.} 
\commHL{By only using some basic statistics of the underlying graphs such as edge density and diameter, our method can achieve perfect accuracy when there is a significant difference in the topology of the underlying graphs between the synchronizing and the non-synchronizing examples. However, in the problem setting where these graph statistics cannot distinguish the two classes very well (e.g., when the graphs are generated from the same random graph model), we find that pairing a few iterations of the initial dynamics along with the graph statistics as the input to our classification algorithms can lead to significant improvement in accuracy; far exceeding what is known by the classical oscillator theory. More surprisingly, we find that in almost all such settings, dropping out the basic graph statistics and training our algorithms with only initial dynamics achieves nearly the same accuracy.} We demonstrate our method on three models of continuous and discrete coupled oscillators --- the Kuramoto model,  Firefly Cellular Automata, and Greenberg-Hastings model. Finally, we also propose an ``ensemble prediction'' algorithm that successfully scales our method to large graphs by training on dynamics observed from multiple random subgraphs. 
\end{abstract}
\begin{document}

\flushbottom
\maketitle
%
%
\thispagestyle{empty}

\section*{Introduction}

Many important phenomena that we would like to understand $-$ formation of public opinion, trending topics on social networks, movement of stock markets, development of cancer cells, the outbreak of epidemics, and collective computation in distributed systems $-$ are closely related to predicting large-scale behaviors in networks of locally interacting dynamic agents. Perhaps the most widely studied and mathematically intractable of such collective behavior is the \textit{synchronization} of coupled oscillators (e.g., blinking fireflies, circadian pacemakers, BZ chemical oscillators), and has been an important subject of research in mathematics and various areas of science for decades \cite{strogatz2000kuramoto, acebron2005kuramoto}. Moreover, it is closely related to the \textit{clock synchronization} problem, which is essential in establishing shared notions of time in distributed systems and has enjoyed fruitful applications in many areas including wildfire monitoring, electric power networks, robotic vehicle networks, large-scale information fusion, and wireless sensor networks
 \cite{dorflfer2012synchronization, nair2007stable, pagliari2010scalable}. 
 
\begin{figure*}[htpb]
	\centering
	\includegraphics[width=1\textwidth]{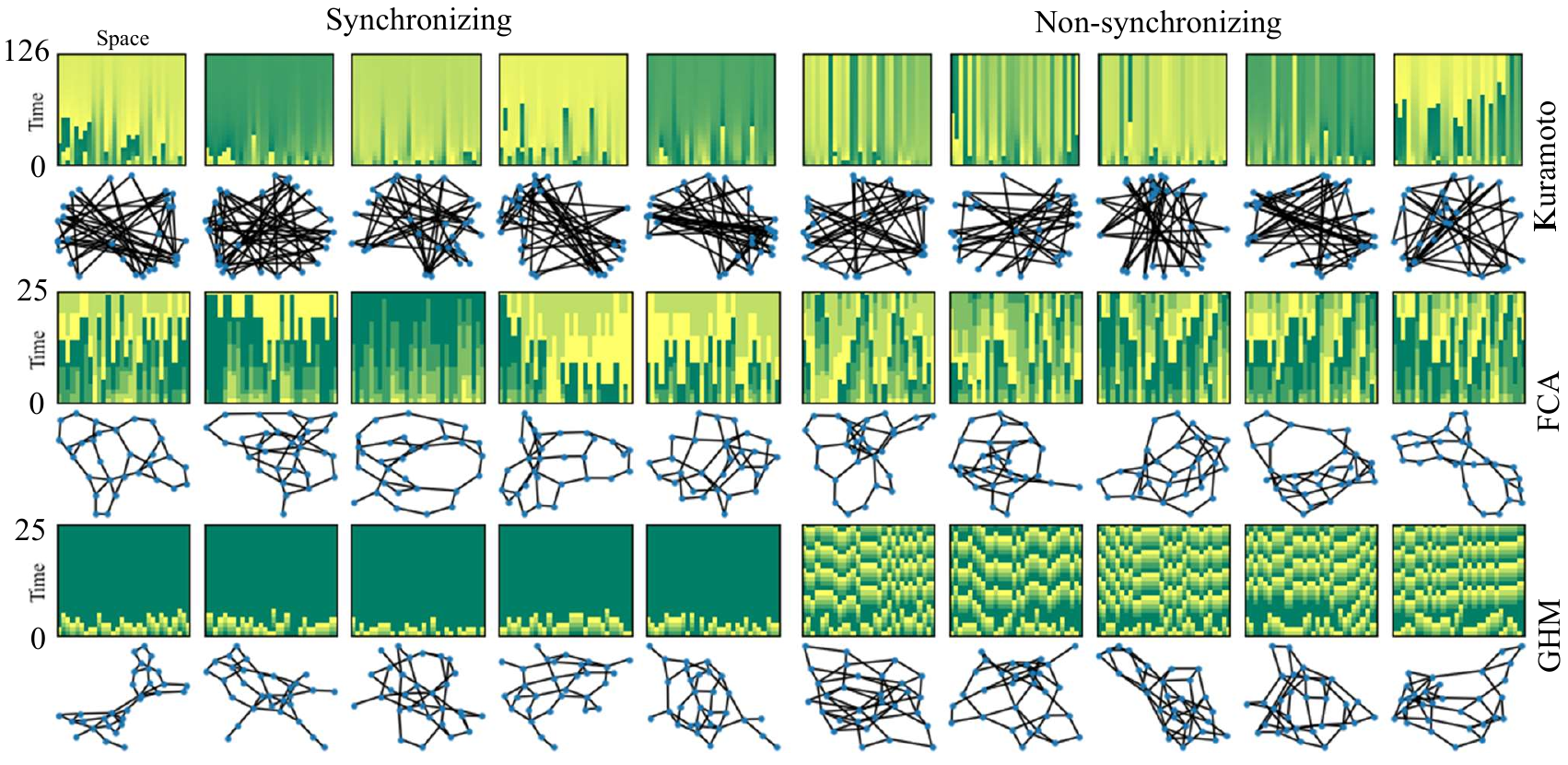}
	\caption{Sample points in the 30-node dynamics dataset for synchronization prediction. The heat maps show phase dynamics on graphs beneath them, where colors represent phases and time is measured by iterations from bottom to top (e.g. $t=0$ to $t=25$). Each example is labeled as `synchronizing' if it synchronizes at iteration 1758 for the Kuramoto model (70 for FCA and GHM) 
	and `non-synchronizing' otherwise. 
	Synchronizing examples have mostly uniform colors in the top row.  For training, only a portion of dynamics is used so that the algorithms rarely see a fully synchronized example (see Figure \ref{fig:dataset_visualization_training}). 
	} 
	\label{fig:dataset_visualization_full}
\end{figure*}

For a system of deterministic coupled oscillators (e.g., the Kuramoto model \cite{kuramoto2003chemical}), the entire forward dynamics (i.e., the evolution of phase configurations) is analytically determined by 1) the initial phase configuration and 2) the graph structure (see Figure \ref{fig:dataset_visualization_full}). In this paper, we are concerned with the fundamental problem of \textit{predicting whether a given system of coupled oscillators will eventually synchronize}, \commHL{using some information on the underlying graph or on the initial dynamics (that is, early-stage of the forward dynamics).} More specifically, we consider the following three types of synchronization prediction problems (see Figure \ref{fig:dataset_visualization_training}): 
\begin{customquestion}{Q1}\label{Q1}
Given the initial dynamics and graph structure, can we predict whether the system will eventually synchronize?
\end{customquestion}

\begin{customquestion}{Q2}\label{Q2}
Given the initial dynamics and not knowing the graph structure, can we predict whether the system will eventually synchronize?
\end{customquestion}

\begin{customquestion}{Q3}\label{Q3}
Given the initial dynamics partially observed on a subset of nodes and possibly not knowing the graph structure, can we predict whether the whole system will eventually synchronize?
\end{customquestion}

Analytical characterization of synchronization would lead to a perfect algorithm for the synchronization prediction problems above. However, while a number of sufficient conditions \commHL{on graph topology \cite{eom2016concurrent,boccaletti2006complex,chowdhury2020effect, lyu2015synchronization} (e.g., complete graphs or trees)}, model parameters (e.g., large coupling strength) or on initial configuration (e.g., phase concentration into open half-circle) for synchronization are known, obtaining an analytic or asymptotic solution to the prediction question, in general, appears to be out of reach, \commHL{especially when these sufficient conditions for synchronization are not satisfied. Namely, we are interested in predicting the synchronization of coupled oscillators where the underlying graph\commHL{s are non-isomorphic} and the initial configuration is not confined within an open half-circle in the cyclic phase space.} 
Since the global behavior of coupled oscillators is built on non-linear local interactions, as the number of nodes increase and the topology of the underlying graphs become more diverse, the behavior of the system becomes rapidly intractable. To provide a sense of the complexity of the problem, note that there are more than $10^{9}$ non-isomorphic connected simple graphs with $11$ nodes \cite{combinatorial_data}. 

However, the lack of a general analytical solution does not necessarily preclude the possibility of successful prediction of synchronization. In this work, we propose a radically different approach to this problem that we call \textit{Learning To Predict Synchronization} (L2PSync), where we view the synchronization prediction problem as a binary classification problem for two classes of `synchronizing' and `non-synchronizing', \commHL{based on the fact that any given deterministic coupled oscillator system eventually synchronizes or converges to a non-synchronizing limit cycle.} \commHL{In this work, we consider three models of continuous and discrete coupled oscillators --- the Kuramoto model (KM) \cite{acebron2005kuramoto}, Firefly Cellular Automata (FCA) \cite{lyu2015synchronization}, and Greenberg-Hastings model (GHM) \cite{greenberg1978spatial}.}

\commHL{Utilizing a few basic statistics of the underlying graphs, our method can achieve perfect accuracy when there is a significant difference in the topology of the underlying graphs between the synchronizing and the non-synchronizing examples (see Figures \ref{fig:toy_graph_dynamics} and \ref{fig:toy}). We find that when these graph statistics cannot separate the two classes of synchronizing and non-synchronizing very well (e.g., when the graphs are generated from the same random graph model, see Tables \ref{table:datasets} and \ref{table:datasets_big}), pairing a few iterations of phase configurations in the initial dynamics along with the graph statistics as the input to the classification algorithms can lead to significant improvement in accuracy. Our methods far surpass what is known by the half-circle concentration principle in classical oscillator theory (see our Methods section). We also find that in almost all such settings, dropping out the basic graph statistics and training our algorithms with only initial dynamics achieves nearly the same accuracy as with the graph statistics. Furthermore, we find that our methods are robust under using incomplete initial dynamics only observed on a few small subgraphs of large underlying graphs. }




\vspace{-0.09cm}
\subsection*{Problem statement}

A graph $G=(V,E)$ consists of sets $V$ of nodes and $E$ of edges. Let $\Omega$ denote the \textit{phase space} of each node, which may be taken to be the circle $\mathbb{R}/2\pi \mathbb{Z}$ for continuous-state oscillators or the color wheel $\mathbb{Z}/\kappa \mathbb{Z}$, $\kappa\in \mathbb{N}$ for discrete-state oscillators. We call a map $X:V\rightarrow \Omega$ a \textit{phase configuration}, and say it is \textit{synchronized} if it takes a constant value across nodes (i.e., $X(v)=Const.$ for all $v\in V$). A \textit{coupling} is a function $\mathcal{F}$ that maps each pair $(G,X_{0})$ of graph and initial configuration $X_{0}:V\rightarrow \Omega$ deterministically to a \textit{trajectory} $(X_{t})_{t\ge 0}$ of phase configurations $X_{t}:V\rightarrow \Omega$. For instance, $\mathcal{F}$ could be the time evolution rule for the KM, FCA, or GHM.  Throughout the paper, $\mathbf{1}(\cdot)$ denotes the indicator function.  The main problem we investigate in this work is stated below: 
\begin{description}
\item[$\bullet$] (Synchronization Prediction)
    \textit{Fix parameters $n\in \mathbb{N}$, $T\gg r>0$,  and coupling $\mathcal{F}$. Predict the following indicator function $\mathbf{1}(\text{$X_{T}$ is synchronized})$ given the initial trajectory $(X_{t})_{0\le t\le r}$ and optionally also with \commHL{statistics of} graph $G$.}
\end{description}

We remark that as $T$ tends to infinity, the indicator in the problem statement converges to the following indicator function $\mathbf{1}(\text{$X_{t}$ is eventually synchronized})$ which aligns more directly with the initial questions \ref{Q1} and \ref{Q2} than the indicator in the above problem statement, but determining whether a given system will never synchronize in amounts to finding a non-synchronizing periodic orbit, which is computationally infeasible in general. 
See Figure \ref{fig:dataset_visualization_training} for an illustration of the synchronization prediction problem.

\subsection*{Related works}
\label{subsection:related_works}

 There are a number of recent works incorporating machine learning methods to investigate problems on coupled oscillators or other related dynamical systems, which we briefly survey in this section and summarize in Table \ref{table:comparison}.
 
 \commHL{In Fan et al.\cite{fan2021anticipating}, the authors are concerned with identifying the critical coupling strength at which a system of coupled oscillators undergo a phase transition into synchronization, where the underlying graph consists of 2-4 nodes with a fixed topology (e.g., a triangle or a star with three leaves). Guth et al. \cite{guth2019machine} use binary classification methods with surrogate optimization (SO) in order to learn optimal parameters and predictors of extreme events. Their work is primarily concerned with learning whether or not intermittent extreme events will occur in various 1D or 2D partial differential equation models. 
 Similarly, Chowdhury et al. \cite{chowdhury2021extreme} utilize a long-short term memory (LSTM) \cite{hochreiter1997long} network to predict whether or not an extreme event will occur on globally coupled mean-field logistic maps on complete graphs. Thiem et al.\cite{Thiem2020emergent} use Feed-forward neural networks (FFNN) \cite{bishop2006pattern} to learn coarse-grained dynamics of Kuramoto oscillators and recover the classical order parameter. Biccari et al.\cite{Biccari2020stochastic} use gradient descent (GD) and the random batch method (RBM) to learn control parameters to enhance the synchronization of Kuramoto oscillators. Slightly less related work is Hefny et al.\cite{hefny2015supervised}, where the authors use hidden Markov models, LASSO regression, and spectral algorithms for learning lower-dimensional state representations of dynamical systems and apply their method to a knowledge tracing model for a dataset of students' responses to a survey. }

\begin{table}[htbp]
		\centering
		\begin{tabular}{c|ccccccc}
			\textit{References} & $\#$ nodes & $\#$ graphs  &\# configs. & model & ML &  Goal \\ 
			\hline
			\rule{0pt}{1.1\normalbaselineskip} Fan et al. \cite{fan2021anticipating} & 2-4 & 1 & 1 & Lorenz, KM & FFNN  & Phase transition \\
			\rule{0pt}{1.1\normalbaselineskip} Guth et al. \cite{guth2019machine} & N/A & N/A & 1 & 1D and 2D PDE &  SO  & Extreme events \\
			\rule{0pt}{1.1\normalbaselineskip} Chowdhury et al.\cite{chowdhury2021extreme} & 200 & 1 & 1 & $\begin{matrix} \textup{Mean-field} \\ \textup{logistic map}  \end{matrix}$ & LSTM & Extreme events \\
		    \rule{0pt}{1.1\normalbaselineskip} Thiem et al.\cite{Thiem2020emergent} & 1500-8000 & 1 & 2000 & KM & FFNN & Parameter estimation\\ 
		    \rule{0pt}{1.1\normalbaselineskip} Biccari et al.\cite{Biccari2020stochastic} & 10-1000 & 1 & 1 & KM & GD, RBM & Parameter estimation\\
		    \rule{0pt}{1.1\normalbaselineskip}  Itabashi et al.\cite{itabashi2021evaluating} & 128-256 & 1  &100 & KM & TDA & $\begin{matrix} \textup{Synchronization}/ \\ \textup{Non-synchronization}/ \\ \textup{Chimera state} \end{matrix} $ \\
		    \hline
		    &&&& KM, & RF, GB, &  \text{Synchronization}/ \\
		    This work & 15-600 & \textbf{2K-200K} & 1 & FCA, & FFNN, & \textup{Non-synchronization}\\
		    &&&& GHM & LRCN \\
		\end{tabular}%
		\caption{Comparison of settings in related works on learning coupled oscillator dynamics using machine learning methods. Recent works in the table focus on learning features of dynamics on fixed graphs. \commHL{In contrast, we aim to classify the long-term dynamics of a given system on a diverse set of underlying graphs.} The column for `\# configs.' refers to the number of distinct initial phase configurations considered for each graph in training.
		}
		\label{table:comparison}
	\end{table}

\commHL{ On the other hand, Itabashi et al.\cite{itabashi2021evaluating} consider classifying coupled Kuramoto oscillators according to their future dynamics from certain features derived (using topological data analysis (TDA)) from their early-stage dynamics (or `initial dynamics' in our terminology). As in the other references above, the underlying graph is fixed in each classification task (a complete graph with weighted edges, which may have two or four communities). But unlike in the references above, the initial configuration, instead of the model parameter, is varied to generate different examples on the same underlying graph. The authors observed that some long-term dynamical properties (e.g., multi-cluster synchronization) can be predicted by the derived features of the initial dynamics. This point is consistent with one of our findings that the first few iterations of the initial dynamics may contain crucial information in predicting the long-term behavior of coupled oscillator systems. }
 
 \commHL{While sharing high-level ideas and approaches with the aforementioned works, our work has multiple distinguishing characteristics in the problem setting and approaches. First, in all of the aforementioned works, a dynamical system on a fixed underlying graph is considered. But in our setting, there are as many as 200K non-isomorphic underlying graphs and we seek to predict whether a given system of oscillators on highly varied or even unknown graphs will eventually synchronize or not. \commHL{Furthermore, we also consider the case when machine learning algorithms are trained on partial observation (e.g., initial dynamics restricted on some subgraphs). }
 Second, only in our work, are discrete models of coupled oscillators (e.g., FCA and GHM) considered, whereas only models with continuous phase space (e.g., the KM) are considered in the literature. Third, only in our work, is the classical concentration principle (a.k.a. the `half-circle concentration', see our Methods section) in the oscillator theory used as a rigorous benchmark to evaluate the efficacy of employed machine learning methods. Finally, we employ various binary classification algorithms such as Random Forest (RF) \cite{breiman2001random} Gradient Boosting (GB)\cite{friedman2002stochastic},  Feed-forward Neural Networks (FFNN) \cite{bishop2006pattern}, and our own adaptation of a Long-term Recurrent Convolutional Networks (LRCN) \cite{donahue2015long} . }

 We remark that there are a number of cases where rigorous results are available for the question of predicting the long-term behavior of coupled oscillators on a graph $G$ and initial configuration $X_{0}$. For instance, the $\kappa=3$ instances of GHM and another related model called Cyclic Cellular Automata (CCA) \cite{fisch1990cyclic} have been completely solved \cite{gravner2018limiting}. Namely, given the pair $(G,X_{0})$, the trajectory $X_{t}$ synchronizes eventually if and only if the discrete vector field on the edges of $G$ induced from $X_{0}$ is conservative (see \cite{gravner2018limiting} for details). Additionally, the behavior of FCA on finite trees is also well-known: given a finite tree $T$ and $\kappa\in \{3,4,5,6\}$, every $\kappa$-color initial configuration on $T$ synchronizes eventually under $\kappa$-color FCA if and only if the maximum degree of $T$ is less than $\kappa$; for $\kappa\ge 7$, this phenomenon does not always hold\cite{lyu2015synchronization, lyu2016phase}. \commHL{This theoretical result on FCA was used in the experiment in Figures \ref{fig:toy_graph_dynamics} and \ref{fig:toy}}. Furthermore, there is a number of works on the clustering behavior of these models on the infinite one-dimensional lattice, $\mathbb{Z}$ (FCA\cite{lyu2015synchronization,lyu2019persistence}, CCA\cite{fisch1990one, fisch1991cyclic, fisch1992clustering,lyu2019persistence} and GHM \cite{lyu2019persistence, durrett1991some}).

\section*{Methods}

The pipeline of our approach is as follows. Namely, 1) fix a model for coupled oscillators; 2) generate a \textit{dynamics dataset} of a large number of non-isomorphic graphs with an even split of synchronizing and non-synchronizing dynamics; 3) train a selected binary classification algorithm on the dynamics dataset to classify each example (initial dynamics with or without underlying features of the graph) into one of two classes, `synchronizing' or `non-synchronizing'; 4) validate the accuracy of the trained algorithms on fresh examples by comparing the predicted behavior of the true long-term dynamics. We use the following classification algorithms: Random Forest (RF) \cite{breiman2001random} Gradient Boosting (GB)\cite{friedman2002stochastic},  Feed-forward Neural Networks (FFNN) \cite{bishop2006pattern}, and our own adaptation of Long-term Recurrent Convolutional Networks (LRCN) \cite{donahue2015long} which we call the \textit{GraphLRCN} (further information such as implementation details and hyperparameters can be found in the SI).

As a baseline for our approach, we use a variant of the well-known ``concentration principle'' in the literature on coupled oscillators. Namely, regardless of the details of graph structure and model, synchronization is guaranteed if the phases of the oscillators are concentrated in a small arc of the phase space \textit{at any given time} (see the next subsection). This principle is applied at each configuration up to the training iteration used to train the binary classifiers.

For question \ref{Q3} on synchronization prediction on the initial dynamics partially observed on subgraphs, as well as reducing the computational cost of our methods for answering \ref{Q1} and \ref{Q2}, we propose an ``ensemble prediction'' algorithm (Algorithm \ref{algorithm:collective_predictor}) that scales up our method to large graphs by training on dynamics observed from multiple random subgraphs. Namely, suppose we are to predict the dynamics of some connected $N$-node graphs, where only the initial dynamics are observed on a few small connected subgraphs of $n \ll N$ nodes. We first train a binary classification algorithm on the dynamics observed from those subgraphs and then aggregate the predictions from each subgraph (e.g., using a majority vote) to get a prediction for the full dynamics on $N$ nodes.

\subsection*{The concentration principle for synchronization and baseline predictor}
In the literature on coupled oscillators, there is a fundamental observation that concentration (e.g., into an open half-circle) of the initial phase of the oscillators leads to synchronization for a wide variety of models on arbitrary connected graphs (see, e.g., \cite[Lem 5.5]{lyu2018global}). This is stated in the following bullet point for the KM and FCA and we call it the ``concentration principle''. This principle has been used pervasively in the literature of clock synchronization \cite{nishimura2011robust, klinglmayr2012guaranteeing, proskurnikov2016synchronization, nunez2014synchronization}  and also in multi-agent consensus problems \cite{moreau2005stability, papachristodoulou2010effects, chazelle2011total}. 
\begin{description}
    \item[$\bullet$] (Concentration Principle) \textit{Let $G$ be an arbitrary connected graph. For the Kuramoto model (see eq. (1) in SI) with identical intrinsic frequency and for FCA (see eq. (3) in SI), given dynamics on $G$ synchronize if all phases at any given time are confined in an open half-circle in the phase space $\Omega$. Furthermore, if all states used in the time-$t$ configuration $X_{t}$ are confined in an open half-circle for any $t\ge 1$, then the trajectory on $G$ eventually synchronizes.}
\end{description}
The `open half-circle' refers to any arc of length $<\pi$ for the continuous \commHL{phase} space $\Omega=\mathbb{R}/2\pi \mathbb{Z}$ and any interval of $<\kappa/2$ consecutive integers ($\textup{mod}\,\, \kappa$) for the discrete \commHL{phase} space $\Omega=\mathbb{Z}/\kappa \mathbb{Z}$. This is a standard fact known to the literature and it follows from the fact that the couplings in the statement monotonically contract given any initial phase configuration under the half-circle condition toward synchronization. It is not hard to see the above half-circle concentration does not hold for GHM. Accordingly, for GHM, we say a phase configuration $X_{t}$ is \textit{concentrated} if $X_{t}$ is synchronized. 

We now introduce the following baseline synchronization predictor: Given $(X_{t})_{0\le t \le r}$ and $T>r$,
\begin{description}
    \item[$\bullet$] (Baseline predictor) \textit{Predict synchronization of $X_{T}$ if $X_{t}$ is concentrated for any $1\le t \le r$. Otherwise, flip a fair coin. }
\end{description}
Notice that the baseline predictor never predicts synchronization incorrectly if $X_{r}$ is concentrated. For non-concentrated cases, the baseline does not assume any knowledge and gives a completely uninformed prediction. Quantitatively, suppose we have a dataset where $\alpha$ proportion of \commHL{entire} samples are synchronizing \commHL{(in all our datastes, $\alpha=0.5$)}. Suppose we apply the baseline predictor where we use the first $r$ iterations of dynamics for each sample. Let $x=x(r)$ denote the proportion of synchronizing samples \commHL{that concentrate by iteration $r$ among all among all synchronizing samples.} Then the baseline predictor's accuracy is given by $x\alpha + (1-\alpha+(1-x)\alpha)/2=0.5+x\alpha/2$, where the term $x\alpha/2$ can be regarded as the gain obtained by using the concentration principle layered onto the uninformed decision. 

\subsection*{An illustrative example: the Kuramoto model on complete vs. ring; GHM on path vs. complete; FCA on tree vs. ring}

We first give a simple example to illustrate our machine learning approach to the synchronization prediction problem. More specifically, we create \commHL{datasets} of an equal number of synchronizing and non-synchronizing examples, where there is a significant difference in the topology of the underlying graphs between the synchronizing and the non-synchronizing examples. In such a setting, one can expect that knowing the basic graph features --- the number of edges, min/max degree, diameter, and the number of nodes --- will be enough to distinguish between synchronizing and non-synchronizing examples.

\begin{figure*}[htpb]
	\centering
	\includegraphics[width=1\textwidth]{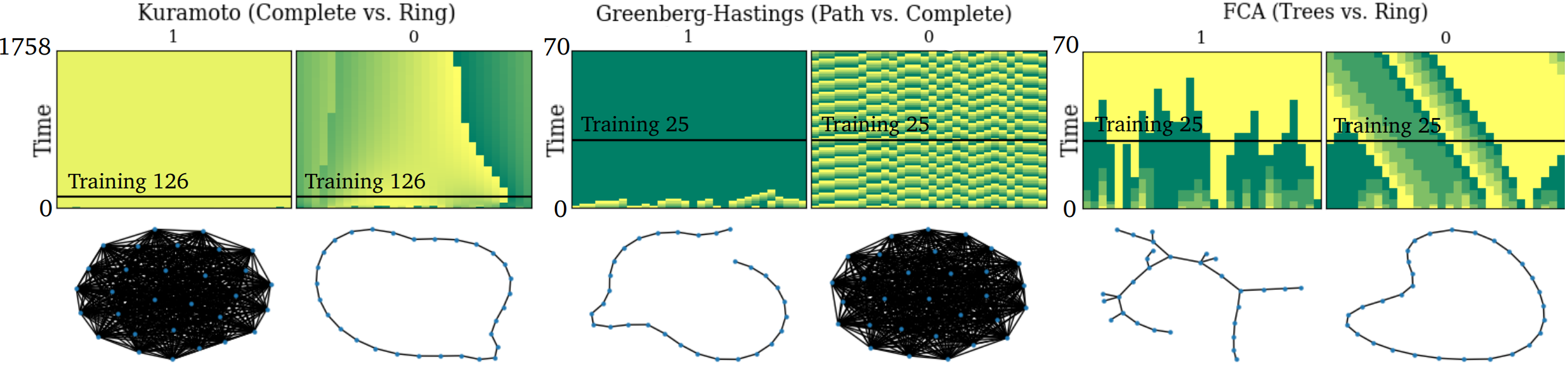}
	\caption{ Examples of synchronizing and non-synchronizing dynamics of Kuramoto, Greenberg-Hastings, and FCA oscillators on 30-node graphs with two topologies (complete vs. ring, path vs. complete, and tree vs. ring, respectively). The heat maps show phase dynamics on graphs beneath them, where colors represent phases and time is measured by iterations from bottom-to-top (e.g. $t=0$ to $t=70$). Synchronizing examples have uniform color in the top row. The horizontal bar indicates the `training iteration', which is the maximum number of iterations in the initial dynamics fed into the classification algorithm for prediction. 
	}
	\label{fig:toy_graph_dynamics}
\end{figure*} 

\begin{figure}[H]
	\centering
	\includegraphics[width=1\textwidth]{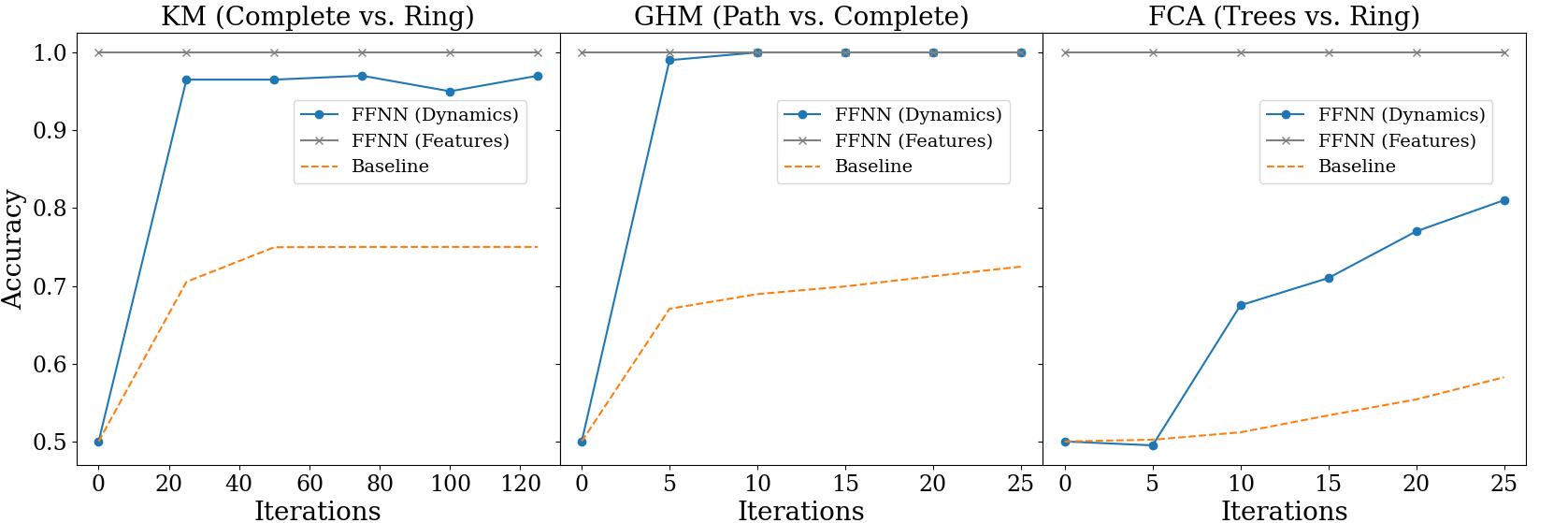}
	\caption{\commHL{Binary classification accuracies of synchronizing vs. non-synchronizing dynamics when the underlying graphs in either of these two classes share the same topology. For the Kuramoto model, we compare complete graphs (synchronizing) to rings (non-synchronizing), for the GHM model we compare rings (synchronizing) vs complete graphs (non-synchronizing), and for FCA we compare trees with maximum degree 4 (synchronizing) to rings (non-synchronizing). FFNN (Features) uses FFNN as the binary classification algorithm with five basic graph features --- the number of edges, min/max degree, diameter, and the number of nodes --- as the input, which achieves perfect classification in all cases. FFNN (Dynamics) uses FFNN as the binary classification algorithm with a varying number of iterations of initial dynamics as input (specified as the horizontal axis), without any information about the underlying graph. The baseline predictor uses the concentration principle (see the main text for more details). Notice that FFNN (Dynamics) far surpasses the baseline.}
	} 
	\label{fig:toy}
\end{figure}

\commHL{For each of the coupled oscillator models of the KM, FCA, and GHM, we create a dataset that consists of 1K synchronizing examples and 1K non-synchronizing examples on 30 node graphs, where all 1K examples in each of the two classes share the same underlying graph topology. For the KM, the synchronizing and the non-synchronizing examples are on a complete graph and on a ring, respectively. For GHM (with $\kappa=5$), the synchronizing and the non-synchronizing examples are on a path and on a complete graph, respectively. Lastly for FCA (with $\kappa=5$), the synchronizing and the non-synchronizing examples are on a tree with a maximum degree at most four and on a ring, respectively. Our choice of these graph topologies reflects rigorously established results in the literature. Namely, it is well-known that coupled Kuramoto oscillators (with identical intrinsic frequencies) on a complete graph will always synchronize, and one can easily generate a non-synchronizing initial configuration for Kuramoto oscillators on a ring \cite{strogatz2000kuramoto,acebron2005kuramoto}. For GHM, in the SI, we prove that a $\kappa$-color (arbitrary $\kappa \geq 3$) GHM on a path of length $n$ will always synchronize by iteration $n+\kappa$, regardless of the initial configuration. Lastly for FCA, it is known that FCA with $\kappa=5$ color dynamics always synchronize on trees with maximum degree at most four \cite{lyu2015synchronization}. The results of our experiments for these three datasets are shown in Figure \ref{fig:toy}.}

\commHL{For the three datasets described above, we achieve perfect classification accuracy with FFNN of distinguishing synchronizing vs non-synchronizing by using the following five basic statistics of the underlying graphs --- 
the number of edges, min/max degree, diameter, and the number of nodes  --- as expected by the theoretical results we mentioned above. Additionally, we hide the graph statistics completely and train FFNN only using the initial dynamics up to a variable training iteration. As we feed in more initial dynamics for training, FFNN quickly improves in prediction accuracy, far more rapidly than the baseline does. This indicates that FFNN may be quickly learning some latent distinguishing features from the initial dynamics that are more effective than the half-circle concentration employed by the baseline predictor.}

\subsection*{Generating the dynamics datasets}


\commHL{In the example we considered in Figures \ref{fig:toy_graph_dynamics} and \ref{fig:toy}, there was a clear distinction between the topology of the underlying graphs that were synchronizing and the graphs that were non-synchronizing, and basic graph statistics can yield perfect classification accuracy through FFNN as the choice of the binary classifier. In this section, we consider datasets for which it is much more difficult to classify based only on the same graph statistics. The way we do so is by generating a large number of underlying graphs from the \textit{same} random graph model with the \textit{same} parameters. In this way, even though individual graphs realized from the random graph model may have different topologies, graph statistics such as edge density or maximum degree are concentrated around their expected values. This is in contrast to the set of underlying graphs in the three datasets we considered in Figures \ref{fig:toy_graph_dynamics} and \ref{fig:toy}, where half of them are isomorphic to each other and another half are also isomorphic to each other. For the datasets we generate in this way, which we will describe in more detail shortly after, classifying only with the basic graph statistics achieves an accuracy of 60-70\% (in comparison to 100\% as before).
}

\commHL{We generate a total of tweleve datasets described in Tables \ref{table:datasets} and \ref{table:datasets_big} as follows.} Data points in each dataset consist of three statistics computed for a pair $(G,X_{0})$ of an underlying graph, $G=(V,E)$, and initial configuration, $X_{0}:V\rightarrow \Omega$: 1) first $r$ iterations of dynamics $(X_{t})_{0\le t \le r}$ (using either the KM, FCA, or GHM), (optional) 2) features of $G$ and $X_{0}$, and 3) the label that indicates whether $X_{T}$ is concentrated or not. We optionally include the following six features: \textit{number of edges}, \textit{min degree}, \textit{max degree}, \textit{diameter}, \textit{number of nodes}, and \textit{quartiles of initial phases in $X_{0}$}. We say a data point is `synchronizing' if the label is 1, and `non-synchronizing,' label 0, otherwise. Every dataset we generate contains an equal number of synchronizing and non-synchronizing examples, and the underlying graphs are all connected and non-isomorphic \commHL{(as opposed to the datasets in Figures \ref{fig:toy_graph_dynamics} and \ref{fig:toy})}.



    To generate a single $n$-node graph, we use an instance of the Newman-Watts-Strogatz (NWS) model\cite{newman2002random}, which originally has three 
    parameters $n$ (number of nodes), $p$ (shortcut edge probability) and $k$ (initial degree of nodes; We use the implementation in \texttt{networkx} python package, see \cite{hagberg2008exploring}), with an added integer parameter $M$ (number of calls for adding shortcut edges). Namely, we start from a cycle of $n$ nodes, where each node is connected to its $k$ nearest neighbors. Then we attempt to add a new edge between each initial non-edge $(u,v)$ with probability $p/(n-k-3)$ independently $M$ times. The number of new edges added in this process follows the binomial distribution with $\binom{n}{2} - \frac{nk}{2}$ trials with success probability $\left(1 - \left(1-\frac{p}{n-k-1} \right) \right)^{M}\approx pm/(n-k-1)$. This easily yields that the expected number of edges in our random graph model is  $\frac{nk}{2}  +  \frac{n^{2}pM}{2(n-k-1)} + O(k)$.

\begin{table}[htbp]
	\centering
	\begin{tabular}{ccccccc}
		\hline 
		 \textit{Datasets} & $\dataset{KM}_{15}$ & $\dataset{KM}_{30}$ & $\dataset{FCA}_{15}$ & $\dataset{FCA}_{30}$ & $\dataset{GHM}_{15}$ & $\dataset{GHM}_{30}$  \\
		\hline 
		\# nodes & 15 & 30 & 15 & 30 &15 & 30 \\
		avg of \# edges & 29.65 &57.49 & 23.91 & 47.45 & 22.88& 48.18 \\
		std of \# edges & 3.42 & 5.67 & 2.34& 4.15 & 2.35& 4.11 \\
		avg diameter & 4.32 & 7.00 & 4.32 & 6.01 & 4.30 & 6.12 \\
		std of diameter & 0.68 & 1.29 & 0.67 & 0.95 & 0.66 & 0.98\\
        $r$ (training iter) & 126 & 126 & 25 & 25 & 25 & 25 \\
        $T$ (prediction iter) & 1758 &  1758 & 70 & 70 & 70 & 70 \\
        \# Sync.   & 100K & 40K & 100K & 40K & 100K & 40K\\
        \# Nonsync.  & 100K & 40K & 100K & 40K & 100K & 40K\\
        \hline
	\end{tabular}%
	\caption  { 
Dynamics datasets were generated for three models with two node counts. In each dataset, all graphs are connected and non-isomorphic. $\#$ Sync. denotes the number of examples in the dataset such that the phase configuration $X_{T}$ at iteration $T$ is concentrated}
	\label{table:datasets}
\end{table}

\begin{figure*}[htpb]
	\centering
	\includegraphics[width=1\textwidth]{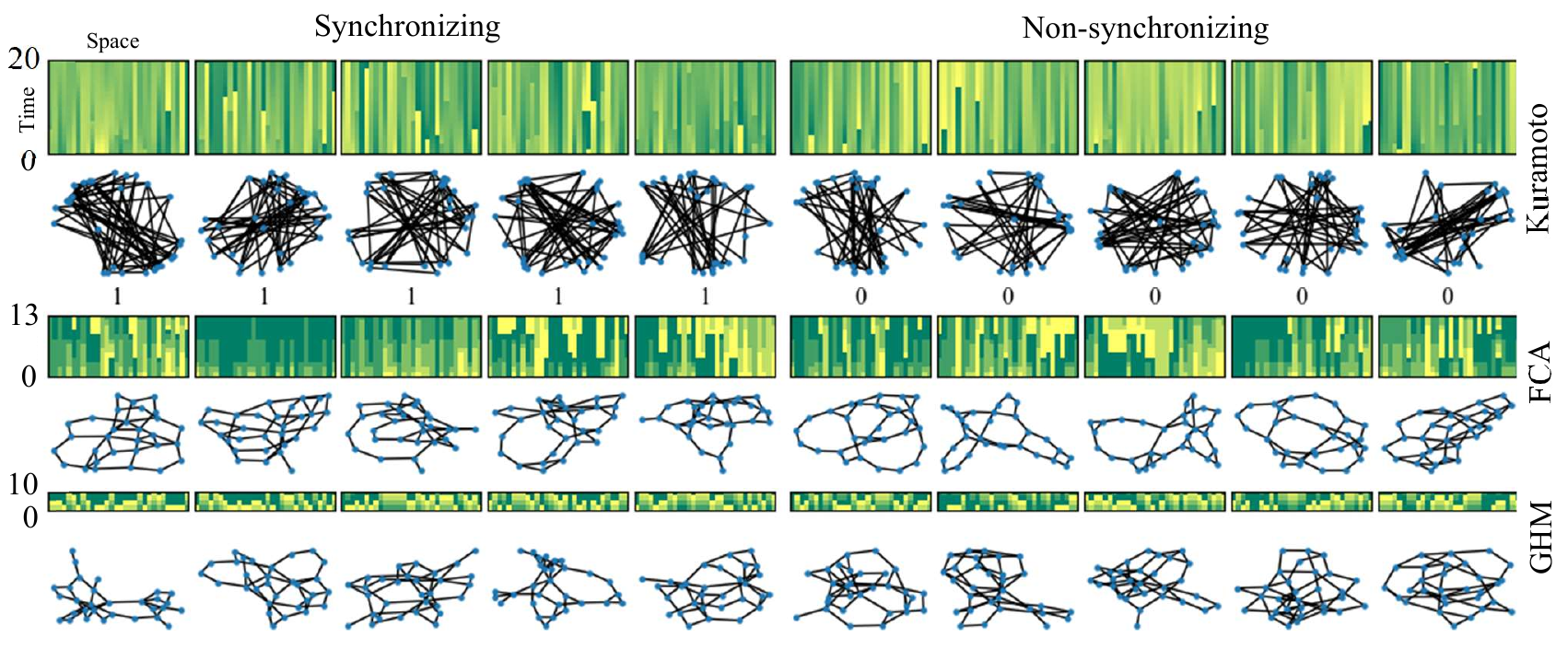}
	\caption{Sample points in the 30-node training data set for synchronization prediction. The full dataset consists of 40K synchronizing and 40K non-synchronizing 30-node connected non-isomorphic graphs and dynamics on them for each of the three models the KM, FCA, and GHM (\commHL{see} Table \ref{table:datasets}). 
	} 
	\label{fig:dataset_visualization_training}
\end{figure*} 

    The NWS model of random graphs is known to exhibit the `small world property', which is to have relatively small mean path lengths in contrast to having a low local clustering coefficient. This is a widely observed property of many real-world complex networks \cite{newman2002random} --- as opposed to the commonly studied Erd\"os-R\'eyni graphs. It is known that many models of coupled oscillators have a hard time synchronizing when the underlying graph is a ring, as the discrepancy between oscillators tends to form traveling waves circulating on the ring \cite{lee2008experiments}. On the other hand, it is observed that coupled oscillator systems on dense graphs are relatively easier to synchronize. For instance, Kassabov, Strogatz, and Townsend recently showed that Kuramoto oscillators with an identical natural frequency on a connected graph where each node is connected to at least 3/4 of all nodes are globally synchronizing for almost all initial configurations \cite{kassabov2021sufficiently}. Since we intend to generate both synchronizing and non-synchronizing examples to form a balanced dataset, it is natural for us to use a random graph model that sits somewhere between rings and dense graphs. In this sense, NWS is a natural choice for a random graph model for our purpose. Using other models such as Erd\"os-R\'eyni, for example, \commHL{for} generating the balanced dataset of synchronizing and non-synchronizing examples as in Tables \ref{table:datasets} and \ref{table:datasets_big}, \commHL{is} computationally very demanding.

    In Table \ref{table:datasets}, we make note of the average graph edge counts and standard deviations of the random graphs that were used to simulate our models. These characteristics, average edge count and standard deviation of edges counts, correspond to the clustering of edges in the graph, which intuitively affects the propagation of information or states in our cellular automata.


Also in Table \ref{table:datasets}, we give a summary of the six datasets on the three models for two node counts $n=15,30$, each with 200K and 80K examples, respectively, which we refer to as $\dataset{KM}_{n}$, $\dataset{FCA}_{n}$, $\dataset{GHM}_{n}$ for $n=15,30$. Underlying graphs are sampled from the NWS model with parameters $n\in \{15,30\}$, $M=1$, and $p=0.85$ for the KM and $p=0.65$ for FCA and GHM. In all cases, we generated about 400K examples and subsampled 200K and 80K examples for $n=15,30$, respectively, so that there are an equal number of synchronizing and non-synchronizing examples, with all underlying graphs as non-isomorphic. The limits for both sets were chosen by memory constraints imposed by the algorithms used. To give a glance at the datasets, we provide visual representations. In Figure \ref{fig:dataset_visualization_training}, we show five synchronizing and non-synchronizing examples in $\dataset{KM}_{30}$, $\dataset{FCA}_{30}$, and $\dataset{GHM}_{30}$. \\

\begin{table}[htbp]
	\centering
	\begin{tabular}{ccccccc}
		\hline 
		 \textit{Datasets} & $\dataset{FCA}_{600} $& $\dataset{FCA}_{600}'$ & $\dataset{FCA}_{600} ''$ & $\dataset{KM}_{600}$ & $\dataset{KM}_{600}'$ & $\dataset{KM}_{600} ''$\\
		\hline 
		\# nodes & 600 & 600 & 300-600 &  600 & 600 & 300-600\\
		std of \# nodes & 0 & 0 & 86.60 &  0 & 0 & 86.60 \\
		avg of \# edges &2985.53&4749.24&2799.49 &1051.08&1109.85&757.89\\
		std of \# edges &37.85&2371.72& 1461.08&19.66&79.47&56.81\\
        $r$  & 50 & 50 & 50& 400 & 400 & 400 \\
        $T$ & 600 &  600 & 600 & 1758 & 1758 & 1758\\
        \# Sync.   & 1K & 1K & 1K & 1K & 1K & 1K \\
        \# Nonsync.  & 1K & 1K & 1K & 1K & 1K & 1K \\
        \hline
	\end{tabular}%
	\caption{Dynamics datasets generated for FCA and Kuramoto on 600 nodes ($\dataset{FCA}_{600}$, $\dataset{KM}_{600}$, $\dataset{FCA}_{600}'$, $\dataset{KM}_{600}'$) and on 300-600 nodes ($\dataset{FCA}_{600}''$, $\dataset{KM}_{600}''$). In each dataset, all graphs are connected and non-isomorphic. $\#$ Sync. denotes the number of examples in the dataset such that the phase configuration $X_{T}$ at iteration $T$ is concentrated. Here $r$ and $T$ refers to the training and the prediction iterations defined in Table \ref{table:datasets}.} 
	\label{table:datasets_big}
\end{table}

We also generated six dynamics datasets with a larger number of nodes on FCA and Kuramoto dynamics, as described in Table \ref{table:datasets_big}. The fixed node datasets $\dataset{FCA}_{600}$, $\dataset{FCA}_{600}'$, $\dataset{KM}_{600}$, and $\dataset{KM}_{600}'$ each consist of 1K synchronizing and non-synchronizing examples of FCA and Kuramoto dynamics on non-isomorphic graphs of 600 nodes. The underlying graphs for the fixed node FCA datasets are generated by the NWS model with parameters $n=600$, $p=0.6$ and $N=5$ for $\dataset{FCA}_{600}$ and $p\sim\textup{Normal(}\mu_X,0.04)$, $\mu_X\sim\textup{Uniform}(0.32,0.62)$, for each $N\sim \textup{Uniform}(\{1,2,\dots,20\})$ calls for $\dataset{FCA}_{600}'$. Similarly, to generate the fixed node Kuramoto datasets, $\dataset{KM}_{600}$ and $\dataset{KM}_{600}'$, we used the NWS model with parameters $n=600$, $p=0.15$ and $N=5$ for $\dataset{KM}_{600}$ and $p\sim\textup{Normal(}\mu_X,0.04)$, $\mu_X\sim\textup{Uniform}(0.32,0.62)$, for each $N\sim \textup{Uniform}(\{1,2,\dots,20\})$ calls for $\dataset{KM}_{600}'$. Consequently, the number of edges in the graphs from $\dataset{KM}_{600}$ and $\dataset{FCA}_{600}$ are sharply concentrated around its mean whereas $\dataset{KM}_{600}'$ and $\dataset{FCA}_{600}'$ have much greater overall variance in the number of edges (see Table \ref{table:datasets_big}). For the varied node datasets, $\dataset{FCA}_{600}''$ and $\dataset{KM}_{600}''$, we kept $p$ distributed in the same way as done for $\dataset{KM}_{600}'$ $\dataset{FCA}_{600}'$, but for each $N\sim \textup{Uniform}(\{1,2,\dots,20\})$ calls of adding shortcut edges, but additionally varied the number of nodes as $n\sim \textup{Uniform}(\{300,301,\dots,600\})$. In this case, both the number of nodes and edges have relatively greater variation compared to the other datasets.

We omit the GHM from this experiment because the dynamics are extremely prone to non-synchronization as a network has more cycles\cite{durrett1993asymptotic}. Hence, for these large graphs, almost all \commHL{GHM dynamics} will be non-synchronizing. \commHL{For instance, for the same set of networks we used to generate the $\dataset{FCA}_{600}$ dataset (consisting of 40K networks of 600 nodes), none of the GHM dynamics synchronized. Hence, there is no meaningful} classification problem to be discussed since the presence of synchronizing graphs that follow this model is extremely sparse.

\subsection*{Scaling up by learning from \commHL{subgraphs}}

In this subsection, we discuss a way to extend our method for the dynamics prediction problem simultaneously in two directions; 1) larger graphs and 2) a variable number of nodes. The idea is to train our dynamics predictor on subsampled dynamics of large graphs (specified as induced subgraphs and induced dynamics), and to combine the local classifiers to make a global prediction.  In the algorithm below, $f(X_{T}):=\mathbf{1}(\textup{$X_{T}$ is concentrated})$, and if $X_{t}$ is a phase configuration on $G=(V,E)$ and $G_{i}=(V_{i},E_{i})$ is a subgraph of $G$, then $X_{t}|_{G_{i}}$ denotes the restriction $v\mapsto X_{t}(v)$ for all $v\in V_{i}$.

\begin{algorithm}[H]
	\caption{Ensemble Prediction of Synchronization}
	\label{algorithm:collective_predictor}
	\begin{algorithmic}[1]
		\State \textbf{Input:} Dynamics dataset on graphs with $\ge N$ nodes; \,\, Test point $(G', (X'_{t})_{0\le t \le r})$; 
		\State \textbf{Parameters:} $n_{0}\le N$ (size of subgraphs),\,\, $k_{\textup{train}}$, $k_{\textup{test}}$ ($\#$ of subgraphs) ,\,\, $\theta$ (prediction threshold)
		\State \textbf{Subsample Dynamics:} 
		\State \quad For each data point  $(G,(X_{t})_{0\le t \le r}, f(X_{T}))$: 
		\State \quad\quad Sample $n_{0}$-node connected subgraphs $G_{1},\dots,G_{k_{\textup{train}}}$ of $G$;
		\State \qquad Form restricted triples $(G_{i}, (X_{t}|_{G_{i}})_{0\le t\le r}, f(X_{T}))$
		\State \textbf{Train Dynamics Predictor:}
		\State \quad Train a binary classifier on the restricted triples;
		\State \textbf{Ensemble Prediction:} 
		\State \quad Sample $n_{0}$-node connected subgraphs $G_{1}',\dots,G_{k_{\textup{test}}}'$ of $G'$; 
		\State \quad $\hat{f}:=$ mean of predictions of $f(X_{T})$ on subdynamics on $G_{i}'$'s
		\State \textbf{Output:} $\mathbf{1}(\hat{f}>\theta)$ 
	\end{algorithmic}
\end{algorithm}
\section*{Results}

Regardless of the three models of coupled oscillators and selected binary classification algorithm, we find that our method used to address \ref{Q1} and \ref{Q2} on average shows at least a 30\% improvement in prediction performance compared to this concentration-prediction baseline for dynamics on 30 node graphs. In other words, our results indicate that the concentration principle applied at each configuration is too conservative in predicting synchronization, and there might be a generalized concentration principle that uses the whole initial dynamics, which our machine learning methods seem to learn implicitly.

Using Algorithm \ref{algorithm:collective_predictor} for \ref{Q3}, we achieve an accuracy score of over 85\% for predicting the commonly studied Kuramoto model on 600 node graphs by only using four 30-node subgraphs where the corresponding baseline gets $55\%$ accuracy. In particular, we observe the baseline with locally observed initial dynamics tends to misclassify non-synchronizing examples as synchronizing, as locally observed dynamics can concentrate in a short time scale while the global dynamics do not.


\subsection*{Synchronization prediction accuracy for 15-30 node graphs}

\begin{figure*}[htpb]
    \centering
     \includegraphics[width=1\linewidth]{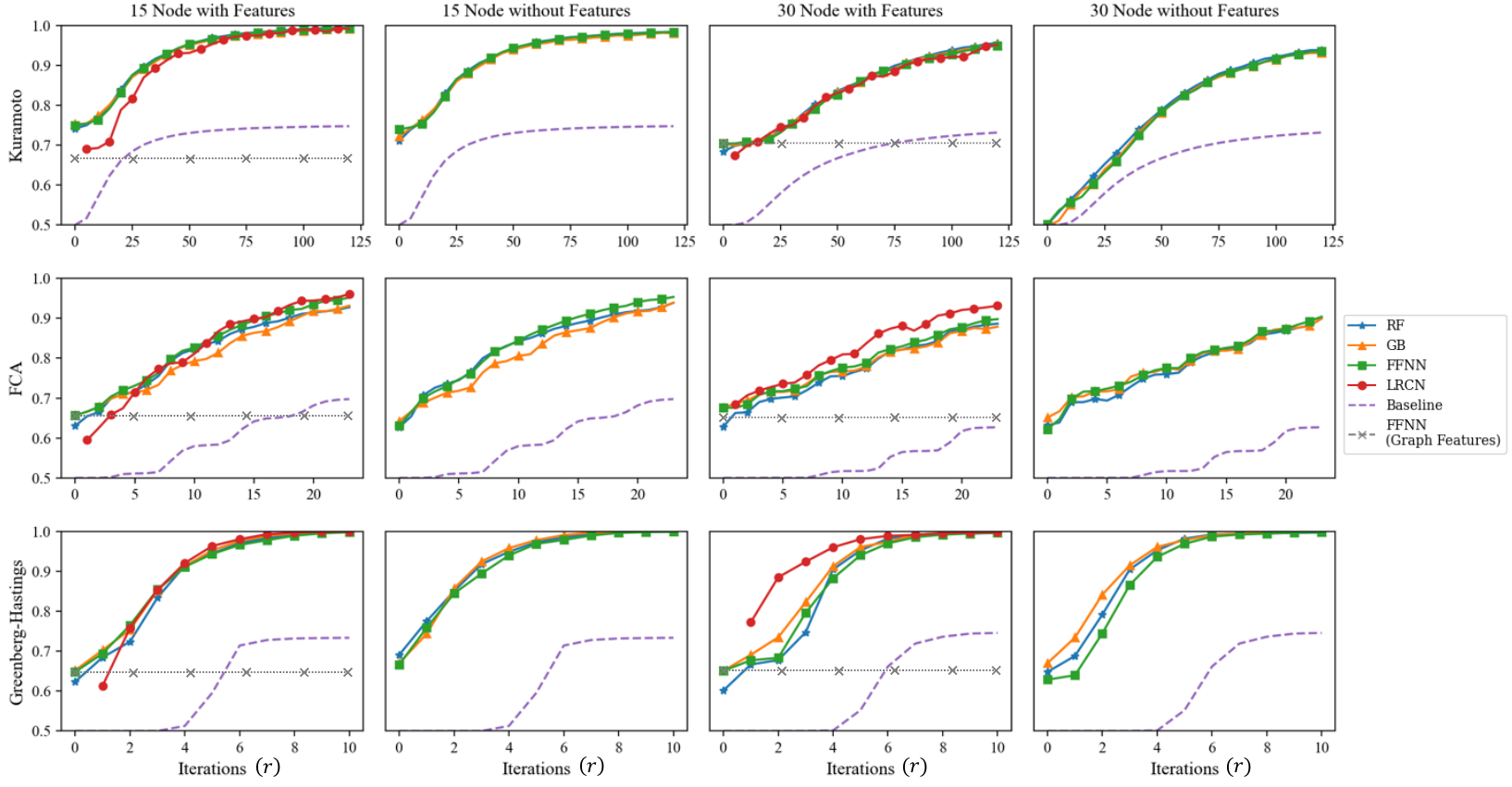}
     \vspace{-0.5cm}
    \caption{Synchronization prediction accuracies of four machine learning algorithms for the KM, FCA, and GHM coupled oscillators synchronization. For each of the six datasets in Table \ref{table:datasets}, we used 5-fold cross-validation with 80/20 split of test/train. Accuracy is measured by the ratio of the number of correctly classified examples to the total number of examples. Algorithms for the second and the fourth columns are trained only with dynamics up to various iterations $r$ indicated on the horizontal axes, whereas the other two columns also use additional graph features. \commHL{$r=0$ in the "with Features" columns indicates the input consisting of initial phase configuration $X_0$ and its \textit{quartiles} paired with the five graph statistics of \textit{number of edges}, \textit{min degree}, \textit{max degree}, \textit{diameter}, \textit{number of nodes}. The gray dashed line represents training FFNN only on the five graph statistics. Hence its accuracy is constant with respect to the varying number of training iterations. } }
    \label{fig:full_plot_accuracy}
\end{figure*}

We apply the four binary classification algorithms for the six 15-30 node datasets described in Table \ref{table:datasets} in order to learn to predict synchronization. Each experiment uses initial dynamics up to a variable number of training iterations $r$ that is significantly less than the test iteration $T$, and the goal is to predict whether each example in the dataset is synchronized at an unseen time $T$. We also experiment with and without the additional graph information described in Table \ref{table:datasets} in order to investigate the main questions \ref{Q1} and \ref{Q2}, respectively. We plot prediction accuracy using four classification algorithms (RF, GB, FFNN, LRCN) and the baseline predictor versus the training iteration $r$, with and without the graph features. The problem of synchronization prediction becomes easier as we increase the training iteration $r$, as indicated by the baseline prediction accuracy in Figure \ref{fig:full_plot_accuracy}. For instance, for $\dataset{KM}_{30}$, there is no trivially synchronizing example at iteration $0$, and but about 10$\%$ of the synchronizing examples will have become phase-concentrated by iteration $r=25$.

Now we discuss the results in Figure \ref{fig:full_plot_accuracy}. \commHL{As we intended before, for the six datasets in Table \ref{table:datasets}, classifying only with the basic graph statistics achieves an accuracy of 60-70\% in comparison to 100\% as in the experiments in Figures \ref{fig:toy_graph_dynamics} and \ref{fig:toy}. This is shown by the gray horizontal lines in Figure \ref{fig:full_plot_accuracy}, which indicates the classification accuracy of FFNN trained only with the same five basic graph statistics used before. This indicates that we may need to use more information on the data points in order to obtain improved classification accuracy. One way to proceed is to use more graph statistics such as clustering coefficients\cite{watts1998collective}, modularity\cite{newman2013spectral}, assortativity \cite{allen2017two}, eigenvalues of the graph Laplacian \cite{zhang2011laplacian}, etc. }

\commHL{Instead, aiming at investigating Question \ref{Q1}, we proceed by additionally using dynamics information, meaning that we include the initial dynamics, $(X_{t})_{0\le t \le r}$, up to a varying number of training iteration $r$.  These results are shown in the first and the third columns in Figure \ref{fig:full_plot_accuracy}. At $r=0$, the input consists of five graph statistics we considered before --- \textit{number of edges}, \textit{min degree}, \textit{max degree}, \textit{diameter}, \textit{number of nodes} --- as well the initial configuration $X_{0}$ with its quartiles. In all cases, all four binary classifiers trained with initial dynamics significantly outperform the \commHL{concentration principle baseline.} The classifiers RF, GB, and FFNN show similar performance in all cases. On the other hand, the GraphLRCN in some cases outperforms the other classifiers, especially so with GHM on 30 nodes. For instance, when $r=20$, $10$ and $4$ for $\dataset{KM}_{30}$, $\dataset{FCA}_{30}$ and $\dataset{GHM}_{30}$, respectively, GraphLRCN achieves a prediction accuracy of $73\%$ (baseline $55\%$, $1.25\%$ concentrates), $84\%$ (baseline $52\%$, $1\%$ concentrates) and $96\%$ (baseline $50\%$, $0\%$ concentrates), respectively. 
}



\commHL{
Now that we have seen that initial dynamics can be used along with the basic graph features to gain a significant improvement in classification performance, we take a step further and see how accurate we can be by only using the initial dynamics, in relation to Question \ref{Q2}. That is, we drop all graph-related features from the input, and train the binary classifiers only with the initial dynamics up to varying iteration $r$. The results are reported in the second and the fourth columns in Figure \ref{fig:full_plot_accuracy}. Since now the classifiers are not given any kind of graph information at all, one might expect a significant drop in the performance. However, surprisingly, except for the dataset $\dataset{KM}_{30}$, the classification accuracies are almost unchanged after dropping graph features altogether from training. 
}

\commHL{
In addition, note that for $\dataset{FCA}_{15}$ and $\dataset{FCA}_{30}$, training all four classifiers on only training iterations of $r=3$ and $r=5$, respectively, without any graph features, produces a prediction accuracy of at least $70\%$. In this case, the baseline achieves only $50\%$, meaning that no synchronizing example is phase-concentrated by a given amount of training iterations. Similarly, for $\dataset{KM}_{15}$, training all four classifiers on only training iterations of $r=5$, without any graph features, produces a prediction accuracy of about $77\%$. In this case, the baseline achieves only $52\%$, so only $10\%$ of all synchronizing examples are phase-concentrated by iteration $r$ (see the formula for baseline accuracy in the section on baseline predictor). This indicates that there may be some evidentiary condition for synchronization on the initial dynamics, which is  \textit{different} from half-circle concentration. A further investigation is due to exactly pin-point what such an evidentiary condition for synchronization might be. }

We give multiple additional remarks on the experiments reported in Figure \ref{fig:full_plot_accuracy}. First, we see that the $\dataset{KM}_{30}$ dataset is adversely affected without the use of the graph features, and beginning at $r=0$, we only achieve 50\% accuracy initially. This is further discussed in the Discussion section along with a systematic analysis of the statistical significance of the features to the classification accuracy in more detail. 

Second, the GraphLRCN binary classifier is offset with respect to the other algorithms as it is fed dynamics information encoded into the adjacency matrix of the underlying graph (see eq. (5) in SI), so it cannot be trained only with the initial dynamics information. Oftentimes, GraphLRCN begins below the classical algorithms but is able to outperform them in intermediate iterations, but by the final training iteration, there is only a negligible difference in the performance over different classifiers.


Third, one may be concerned that the few selected graph features we use for training may not give enough information about the given classification tasks, and whether including further graph features may change the results significantly. We remark that GraphLRCN uses the entire graph adjacency matrix as part of the input (see the SI), so technically it uses every possible graph feature during training. In Figure \ref{fig:full_plot_accuracy}, it does show somewhat improved performance on 30 node graphs with the FCA and the GHM dynamics, but virtually identical results for the Kuramoto dynamics. On the other hand, its performance is diminished for 15 node graphs. This may be due to the fact that the training data for 15 node graphs is not rich enough to properly train GraphLRCN, which is a significantly more complex classification model combining convolutional and recurrent neural networks than the other classification models we use.

Fourth, \commHL{in order to investigate the impact of the actual amount of graphs that were used in the datasets with respect to our experiments}, we perform an additional experiment on the 30 node case, in which we subsample 10K and 40K \commHL{balanced datasets} from our 80K sets to see if we get similar performance to Figure \ref{fig:full_plot_accuracy}, \commHL{where the full datasets of 80K examples are used.} As seen in Figure \ref{fig:volume} in the SI, we can see that most methods \commHL{show} a significantly better accuracy \commHL{than} the baseline, with the exception of the \commHL{GB} and \commHL{FFNN} method\commHL{s} for \dataset{KM}$_{30}$ without graph features. 
We see however that the \commHL{accuracies between different methods become almost identical as we increase the number of training iterations.}
We note that with only a subsampled portion of the datasets being used,  we see larger differences in accuracy between the methods themselves, which contrasts with Figure \ref{fig:full_plot_accuracy}, where all methods were rather saturated and displayed very similar accuracy across all iterations. 

\subsection*{Synchronization prediction accuracy for 300-600 node graphs}

In this section, we investigate Question \ref{Q3}, which is to predict synchronization only based on local information observed on select subgraphs. Note that this is a more difficult task than \ref{Q2}, since not only may we not have information about the underlying graph but we also may not have observed the entire phase configuration. For example, the dynamics may appear to be synchronized at a local scale (e.g., on 30-node connected subgraphs), but there are still large-scale waves being propagated and the global dynamics are not synchronized. Nevertheless, we can use the ensemble prediction method (Algorithm \ref{algorithm:collective_predictor}) to combine decisions based on each subgraph to predict the synchronization of the full graph.

\begin{figure*}[!htpb]
    \centering
     \includegraphics[width=0.85\linewidth]{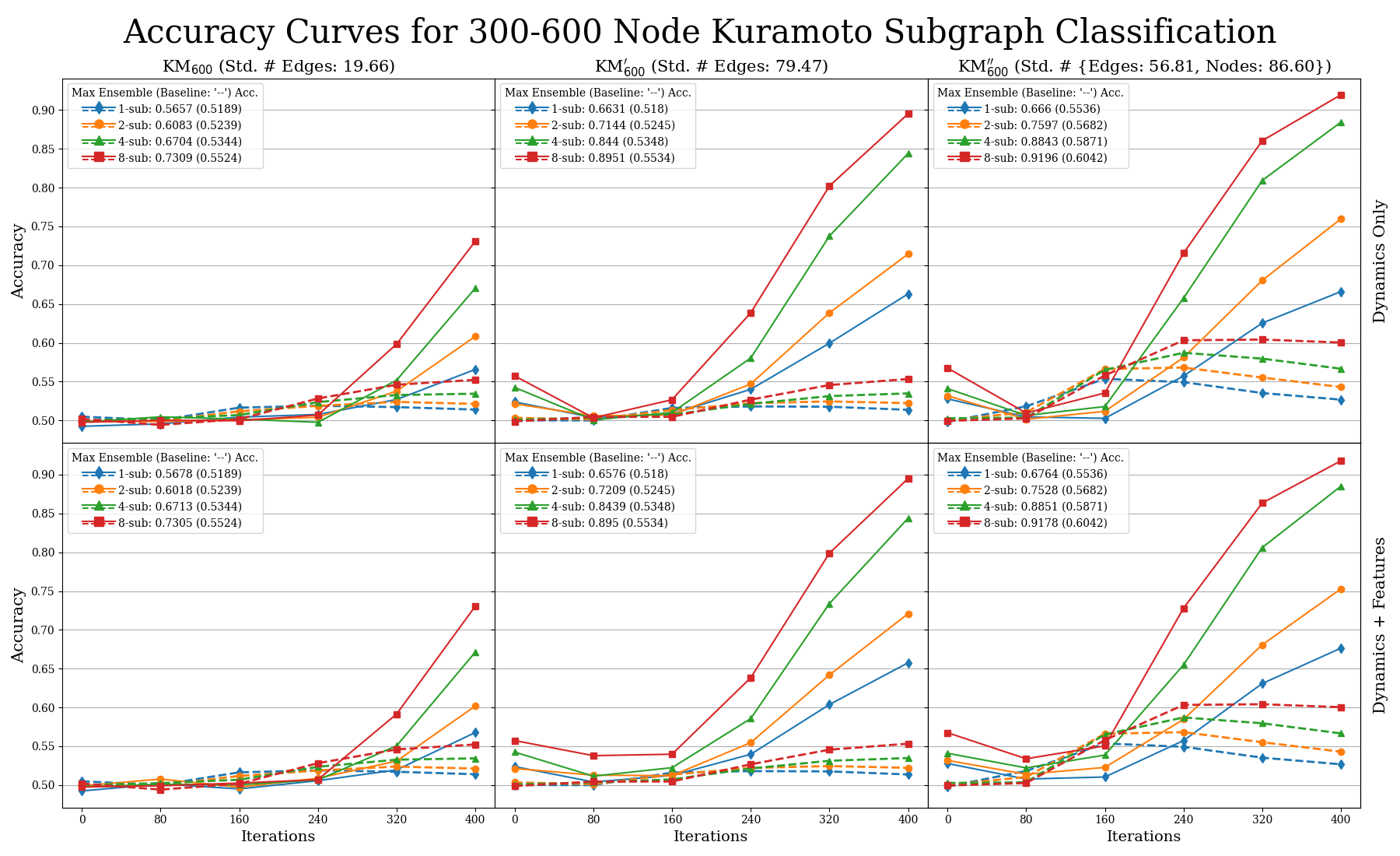}
    \caption{Accuracy curves for predicting synchronization of the Kuramoto model on 600-node graphs from dynamics observed from $k\in\{1,2,4,8\}$ subgraphs of 30 nodes. All plots observe the performance of both the ensemble machine learning (solid) and baseline (dashed) accuracies over increasing amounts of training iterations  $r\in\{0,80,240,320,400\}$. The first row shows results using only dynamics whereas the second row includes both the dynamics and graph features. Maximum accuracies for using $k$ subgraphs are given by `$k$-sub: Acc. (Baseline Acc.)'}
    \label{fig:ensemble_plot_accuracy1}
\end{figure*}

\begin{figure*}[!htpb]
    \centering
     \includegraphics[width=0.85\linewidth]{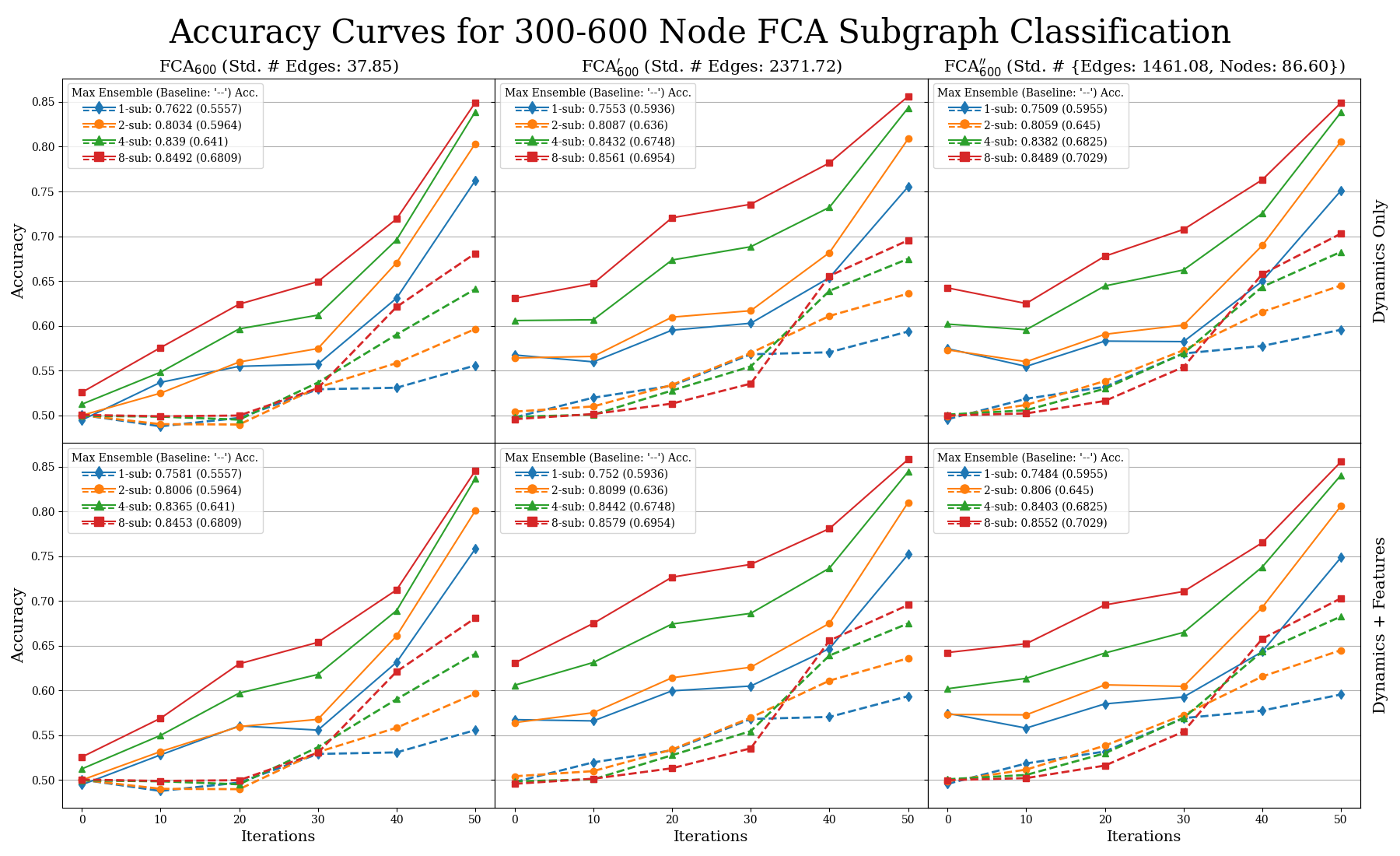}
    \caption{Accuracy curves for predicting synchronization of 5-color FCA on 600-node graphs from dynamics observed from $k\in\{1,2,4,8\}$ subgraphs of 30 nodes. See the caption of Figure \ref{fig:ensemble_plot_accuracy1} for details.}
    \label{fig:ensemble_plot_accuracy2}
\end{figure*}

In Figures \ref{fig:ensemble_plot_accuracy1} and \ref{fig:ensemble_plot_accuracy2}, we report the synchronization prediction accuracy of the ensemble predictor (Algorithm \ref{algorithm:collective_predictor}) on datasets $\dataset{FCA}_{600}$, $\dataset{FCA}_{600}'$, $\dataset{FCA}_{600}''$, and $\dataset{KM}_{600}$, $\dataset{KM}_{600}'$, $\dataset{KM}_{600}''$, respectively, described in Table \ref{table:datasets_big}. We used Algorithm \ref{algorithm:collective_predictor} with $n_{0}=30$ (amount of nodes in the subgraphs) and $k\in \{1,2,4,8\}$ ($\#$ of subgraphs). The binary classification algorithm we used is FFNN. We chose FFNN because as seen in Figure \ref{fig:full_plot_accuracy}, there is no significant difference in accuracies between all methods used in the 30 node case. Furthermore, we can not use the GraphLRCN model, as \commHL{this method relies on knowing} the dynamics and adjacency matrices of the underlying subgraphs. As the baseline, we use a slight modification of the baseline predictor from the 15-30 node classification task. Namely, we combine all phase values observed on all $k$ subgraphs and predict synchronization if they satisfy the concentration principle otherwise we flip a fair coin. Note that the concentration of phases observed on subgraphs does \textit{not} imply synchronization of the full system, as the full phase configuration may not be concentrated.

So far we have been using the metric of accuracy for our synchronization prediction task, which is defined as the ratio of the correctly predicted synchronizing and non-synchronizing examples to the total number of examples in the dataset. There are multiple ways that this metric could be low, for example, if an algorithm is overly conservative and misses lots of synchronizing examples or in the opposite case it may incorrectly classify lots of non-synchronizing examples. \commHL{In order to provide better insights} on these aspects, we also report the performance of our method (Algorithm \ref{algorithm:collective_predictor}) and the baseline in terms of precision and recall. Here, \textit{precision} is formally defined as the proportion of positive classifications that are correct, or $\frac{TP}{TP+FP}$, where $TP$ is "true positive" and $FP$ is "false positive." In our problem, this corresponds to what proportion of predictions in synchronization was truly correct. \textit{Recall} is formally defined as the proportion of accurately identifying true labels in the data, or $\frac{TP}{TP+FN}$, where $FN$ is "false negative." Again with respect to our prediction problem, measuring recall corresponds to how well a given algorithm, ensemble, or baseline, correctly identifies synchronization behavior when presented with synchronizing data.



Our accuracy results for both the ensemble method and baseline are presented in Figures \ref{fig:ensemble_plot_accuracy1} and  \ref{fig:ensemble_plot_accuracy2} for Kuramoto and FCA models, respectively. In these figures, the first columns represents the results for the datasets $\{\dataset{KM,FCA}\}_{600}$ with fixed number of nodes (600) and relatively smaller variation of edge counts (std $\approx 20, \,38$, respectively); the second columns for the datasets $\{\dataset{KM,FCA}\}'_{600}$ with fixed number of nodes (600) and relatively larger variation of edge counts (std $\approx 80,\, 2371$, respectively); and the third for the datasets $\{\dataset{KM,FCA}\}''_{600}$ with varied number of nodes (std $\approx 86$) and relatively larger variation of edge counts (std $\approx 1461, \,56$, respectively). The first rows of both figures represent prediction accuracies using exclusively dynamics data, and the second row utilizes both dynamics and graph features. In the SI, Figures \ref{fig:ensemble_plot_accuracy3} and \ref{fig:ensemble_plot_accuracy4} show the recall and precision curves of the ensemble and baseline methods with the same row and column orientation as the accuracy figures representing different graph datasets and subsets of features. For the Kuramoto data in Figure \ref{fig:ensemble_plot_accuracy1}, we applied the ensemble and baseline algorithms cumulatively up to training iterations $r=0, 80, 160, 240, 320$ and $400$; and for the FCA data in Figure \ref{fig:ensemble_plot_accuracy2}, we applied the ensemble and baseline algorithms cumulatively up to training iterations $r=0, 10, 20, 30, 40$ and $50$ (Note that using $r=0$ means fitting the prediction algorithms at the initial coloring). We additionally remark that these curves were averaged over 30 train-test splits with 80\% training and  20\% testing.

Across all datasets considered, we see that our ensemble method consistently outperforms the baseline method in accuracy at training iterations $r=400$ for Kuramoto and $r=50$ for FCA. For example, across all Kuramoto datasets and feature subsets by the last iteration, $r=400$, the ensemble method on a single subgraph outperforms the baseline algorithm on eight subgraphs; the best that the baseline algorithm does is 60.42\% accuracy compared to the 91.96\% for the ensemble method's best accuracy. For FCA, the baseline accuracy is 70.29\% compared to the best ensemble method score of 85.79\%, both on eight subgraphs. Considering the recall and the precision plots in Figures \ref{fig:ensemble_plot_accuracy3} and \ref{fig:ensemble_plot_accuracy4}  in the SI gives a more detailed explanation. Namely, the ensemble method significantly outperforms the baseline in the recall by at least 35\%, whereas it performs relatively worse in precision than the baseline by at most $10\%$ for both Kuramoto and FCA with $k=8$ subgraphs except for the datasets $\dataset{KM}_{600}$. This means the baseline is `too conservative' in the sense that it misses correctly classifying a large number of synchronizing examples. From this, we deduce that a large number of synchronizing examples exhibit phase concentration over subgraphs much later in the dynamics---making early detection through phase concentration difficult. On the contrary, the high recall scores of the ensemble methods indicate that our method can still detect most of the synchronizing examples by only using local information observed on the selected subgraphs. To elaborate, the baseline only determines phase concentration at a single point in time, whereas the ensemble method is able to learn the whole variation of dynamics up to iteration $r$. Furthermore, between training our data on dynamics only, versus dynamics \textit{and} graph features, across accuracy, recall, and precision curves, the inclusion of graph features hardly improves the maximum score of these values. 

Finally, for both the Kuramoto model and FCA, we observe that high variation of the node and edge counts boosts all performance metrics of accuracy, recall, and precision. For example, for the recall curves we see that for varied edge count, compared to fixed edge count, the recall is considerably higher than fixed edge count by iteration $r=400$, comparing a recall rate of 0.8755 and 0.9759 for 8 subgraphs. Furthermore, the performance gain in introducing a larger variation of node and/or edge counts is significantly larger for the Kuramoto model than for FCA. We speculate that having a larger variation of node and edge counts within the dataset presumably implies a better separation between the synchronizing and the non-synchronizing examples in the space of initial dynamics. We remark that the performance gain here is not due to a better separation between the classes in terms of the graph features, as we can see by comparing the first and the second rows of all figures (Figures \ref{fig:ensemble_plot_accuracy1}, \ref{fig:ensemble_plot_accuracy2} and also Figures \ref{fig:ensemble_plot_accuracy3} and \ref{fig:ensemble_plot_accuracy4}  in the SI).

\section*{Discussion}
\label{subsection:gini}

In Figure \ref{fig:full_plot_accuracy}, observe that not using the additional graph features as input (column 4) decreases the prediction accuracy from $70\%$ to $50\%$ for the case of \dataset{KM}$_{30}$ at initial training iteration $r=0$, but there is no such significant difference for \commHL{the discrete models} \dataset{FCA}$_{30}$ and \dataset{GHM}$_{30}$. In fact, in all experiments we perform in this work (those reported in Figures \ref{fig:full_plot_accuracy},  \ref{fig:ensemble_plot_accuracy1}, \ref{fig:ensemble_plot_accuracy2} and Figures \ref{fig:ensemble_plot_accuracy3} and \ref{fig:ensemble_plot_accuracy4}  in the SI), this is only instance that we see that including the graph features during training affects the prediction accuracy. 

In order to explain why this is the case, we compare the statistical significance of the all features (including the initial coloring) we use for the prediction task. We do so by computing the Gini indexes of all features by repeating the prediction experiment over 300 train/test splits of the datasets \dataset{KM}$_{30}$, \dataset{FCA}$_{30}$, \dataset{GHM}$_{30}$ using GB for the choice of binary classifier. The analysis using GB will be representative of all other binary classifiers since there are negligible differences in their performance in Figure \ref{fig:full_plot_accuracy}. 

\begin{figure}[!hbt]
    \centering
     \includegraphics[width=1\linewidth]{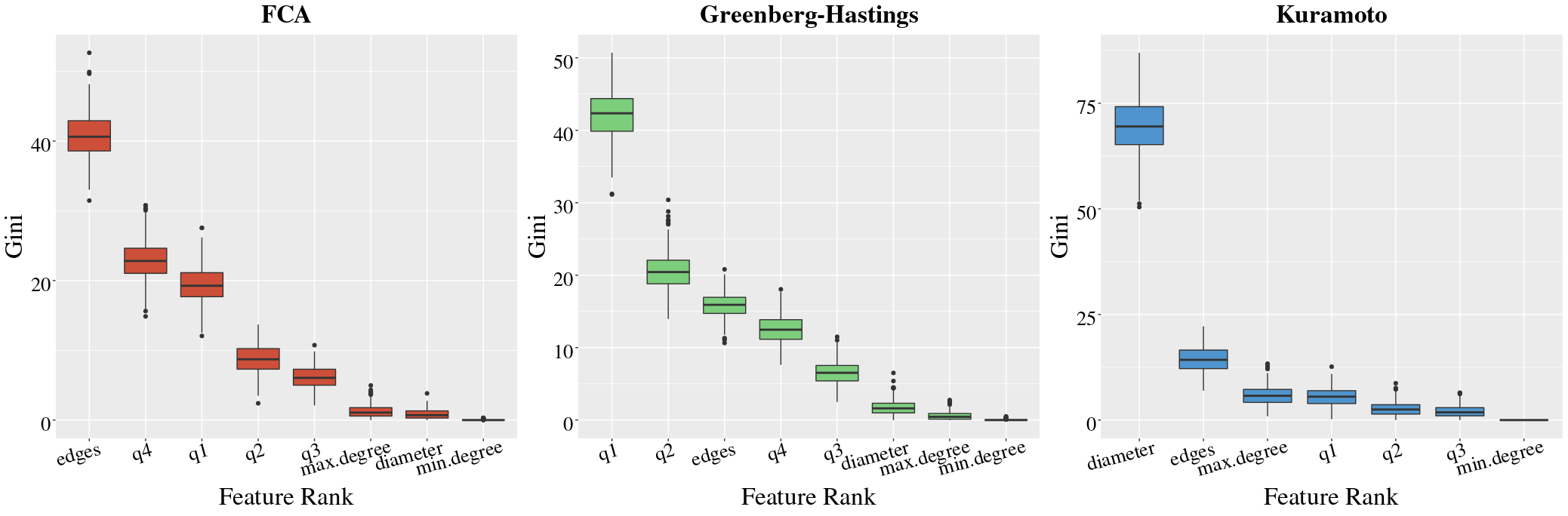}
    \caption{Boxplots of Gini index values sampled from gradient boosting procedure. Color information appears to be very important according to the distribution of Gini index values in both the FCA and GH models, discrete cellular automata models, and diameter overwhelmingly appears to have the greatest importance for Kuramoto.}
    \label{fig:giniboxplot}
\end{figure}

As mentioned before, the procedure works by fitting random subsets of features and iteratively growing decision trees; it also records what is known as the Gini index \cite{hartmann1982application} over each consecutive partition through multiple training iterations. As decision trees split the feature space, the Gini index measures total variance across class labels--- in our case, synchronizing v.s. non-synchronizing --- over each partition. Allowing the supposition that synchronization can be modeled through our graph features data using the gradient boosting method, observing the mean decrease in the Gini index across decision trees allows us to infer feature importance for synchronization over different models and graphs (see Figure \ref{fig:giniboxplot}). See, Hartman et al.\cite{hartmann1982application} for more discussions on GB and computing the Gini index. 


Interestingly, the discrete models of FCA and GHM place significantly greater importance on color statistics such as the initial quartile colorings, and less importance on graph features such as diameter and minimum and maximum degree. Note that the initial color statistics can be directly computed from the initial dynamics even at training iteration $r=0$, so this explains that there is no significant performance difference in prediction accuracy for the discrete models in Figure \ref{fig:full_plot_accuracy}. On the contrary, the Kuramoto model puts greater importance on graph features such as diameter and number of edges rather than the initial color statistics, and such graph features are not available from the information given by the initial dynamics at training iteration $r=0$. From this, we can see why our algorithms show low prediction accuracy in the case of $\dataset{KM}_{30}$. In addition, we note that there is only a negligible difference in prediction accuracy for $\dataset{KM}_{15}$ in Figure \ref{fig:full_plot_accuracy}. We speculate that this is because for $\dataset{KM}_{15}$, as we see in Table \ref{table:datasets}, having only 15 nodes does not allow to have a significant variation in the diameter and the edge counts of the graphs in the dataset. 

Lastly, we remark that the synchronization behavior of the Kuramoto model on complete graphs is well-understood in the literature on dynamical systems 
\cite{strogatz2000kuramoto, kuramoto2003chemical, acebron2005kuramoto}. In our experiments, we observed that graphs with small diameters tend to be synchronized, and this aligns with the literature since graphs with the same number of nodes but with smaller diameters are closer to complete graphs. 

\section*{Conclusion}

\commHL{Predicting whether a given system of coupled oscillators with an underlying arbitrary graph structure will synchronize is a relevant yet analytically intractable problem in a variety of fields. In this work, we offered an alternative approach to this problem by viewing this problem as a binary classification task, where each data point consisting of initial dynamics \commHL{and/or statistics of underlying graphs }needs to be classified into two classes of `synchronizing' and `non-synchronizing' dynamics, depending on whether a given system eventually synchronizes or converges to a non-synchronizing limit cycle. We generated large datasets with non-isomorphic underlying graphs, where classification only using basic graph statistics is challenging. In this setting, we found that pairing a few iterations of the initial dynamics along with the graph statistics as the input to the classification algorithms can lead to significant improvement in accuracy; far exceeding what is known by the half-circle concentration principle in classical oscillator theory. More surprisingly, we found that in almost all such settings, dropping out the basic graph statistics and training our algorithms with only initial dynamics achieves nearly the same accuracy. Finally, we have also shown that our methods are scale well to large underlying graphs by using incomplete initial dynamics only observed on a few small subgraphs . }

Drawing conclusions from our machine learning approaches to the synchronization prediction problem, we pose the following hypotheses:
\begin{description}
    \item[$\bullet$] The entropy of the dynamics of coupled oscillators may decay rapidly in the initial period to the point that the uncertainty of the future behavior from an unknown graph structure becomes not significant.
    
    \item[$\bullet$] The concentration principle applied at any given time is too conservative in predicting synchronization, and there might be a generalized concentration principle that uses the whole initial dynamics, which our machine learning methods seem to learn implicitly.
\end{description}

Given that our machine learning approach is able to achieve high prediction accuracy, we suspect that there may be some analytically tractable characterizations on graphs paired with corresponding initial dynamics signaling eventual synchronization or not, which we are yet to establish rigorously. As mentioned at the end of the Related Works section, previously known characterizing conditions include the initial vector field on the edges induced by the initial color differential 
for the 3-color GHM and CCA \cite{gravner2018limiting}, as well as the number of available states being strictly less than the maximum degree of underlying trees for FCA \cite{lyu2015synchronization, lyu2016phase}. Designing similar target features into datasets and training binary classification algorithms could guide the further analytic discovery of such conditions for the coupled oscillator models considered in this work. 

Furthermore, even though we have focused on predicting only two classes of the long-term behavior of complex dynamical systems as only synchronizing and non-synchronizing dynamics, our method can readily be extended to predicting an arbitrary number of classes of long-term behaviors. For instance, one can consider the $\kappa$-state voter model on graphs, where the interest would be the final dominating color. In such circumstances, one can train $\kappa$-state classification machine learning algorithms on datasets of \commHL{non-isomorphic} graphs. \commHL{We can also consider extending our method to be able to predict synchronization on a network based on parameter control. For instance, we can train many different trajectories on a singular graph using different intrinsic frequencies in the Kuramoto model, and learn to predict what range of values of intrinsic frequencies promote synchronization.} 


Finally, a more ambitious task beyond long-term dynamic behavior quantified by a single metric is the potential extension of our methods to full time-series and graph state regression. In other words, if each node in the graph represents an individual in an arbitrary social network, can we predict the sentiment level for a given topic at any given time $t$ for every single individual in that particular social network? One can again generate large overarching social networks and run many simulations of sentiment dynamics with many possible edge configurations between individuals (for example, measured by the number of mutual friends or likes/shares of posts on social media). The ultimate goal would be a framework for learning to predict, with precision, entire trajectories of complex dynamical systems.

\section*{Acknowledgments}

This work is supported by NSF grant \href{https://www.nsf.gov/awardsearch/showAward?AWD_ID=2010035&HistoricalAwards=false}{DMS-2010035}. We are also grateful for partial support from and the Department of Mathematics and the Physical Sciences Division at UCLA. We also thank Andrea Bertozzi and Deanna Needell for support and helpful discussions. JV is partially supported by Deanna Needell's grant NSF BIGDATA DMS \#1740325.

\section*{Data Availability}

The codes for the main algorithm used during the current study are available in the repository  \url{https://github.com/richpaulyim/L2PSync}. The datasets used and/or analyzed during the current study are available from the corresponding author on reasonable request.


\small{
	\bibliographystyle{amsalpha}
	\bibliography{mybib}
}

\section*{Author contributions statement}
H.B. and R.P.Y. ran simulations and experiments. H.L. led the experiment design. All authors reviewed the manuscript.

\newpage 

\section*{Supplementary Information}

\subsection*{Three models for coupled oscillators}

The \textit{Kuramoto model} (KM) is one of the most well-studied models for coupled oscillators with continuous states  \cite{kuramoto2003chemical, acebron2005kuramoto, strogatz2000kuramoto, lee2010vortices}. Namely, fix a graph $G=(V,E)$ and the continuous state space $\Omega=\mathbb{R}/2\pi\mathbb{Z}$. A given initial phase configuration $X_{0}:V\rightarrow \Omega$ evolves via the following systems of ordinary differential equations 
\begin{align}\label{eq:def_kuramoto}
\frac{d}{dt} X_{t}(v) = \omega_{v} + \sum_{u\in N(v)} K \, \sin(X_{t}(u) - X_{t}(v)),
\end{align}
for all nodes $v\in V$, where $N(v)$ denotes the set of neighbors of $v$ in $G$, $K$ denotes the \textit{coupling strength}, and $\omega_{v}$ denotes the \textit{intrinsic frequency} of $v$. Since we are interested in the dichotomy between synchronization and non-synchronization, we will be assuming identical intrinsic frequencies, which can be assumed to be zero without loss of generality by using a rotating frame. Note that synchronization is an \textit{absorbing state} under this assumption, that is, $X_{t}$ is synchronized for all $t\ge s$ if $X_{s}$ is synchronized. In order to simulate KM, we use the following discretization
\begin{align}\label{eq:def_kuramoto_discrete}
X_{t+h}(v) - X_{t}(v) = h\left(\sum_{u\in N(v)} K \, \sin(X_{t}(u) - X_{t}(v))\right),
\end{align}
with step size $h=0.05$. Accordingly, an `iteration' for KM is a single application of the difference equation \eqref{eq:def_kuramoto_discrete}. After a change of time scale, we write $X_{k}$ for the configuration obtained after $k$ iterations. 

\begin{figure}[H]
	\centering
	\vspace{-0.3cm}
	\includegraphics[width=0.5\textwidth]{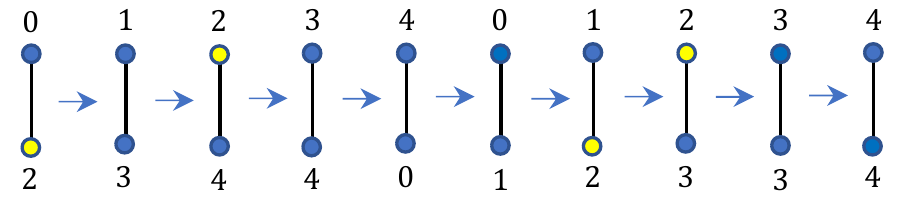}
	\vspace{-0.1cm}
	\caption{An example of 5-color FCA dynamics on two connected nodes. $b(5)=2$ is the blinking color shown in yellow.} 
	\label{fig:FCA_ex1}
\end{figure}

We also consider two discrete models for coupled oscillators. Let $\Omega=\mathbb{N}/\kappa\mathbb{Z}$ be the $\kappa$-state color wheel for some integer $\kappa\ge 3$. The $\kappa$-color \textit{Firefly Cellular Automata} (FCA) is a discrete-state discrete-time model for inhibitory Pulse Coupled Oscillators (PCOs) introduced by Lyu \cite{lyu2015synchronization}. The intuition is that we view each node as an identical oscillator (firefly) of $\kappa$-states, which \textit{blinks} whenever it returns to a designated \textit{blinking color} $b(\kappa)=\lfloor \frac{\kappa-1}{2} \rfloor$; nodes with post-blinking color $c\in (b(\kappa),\kappa-1]$ in contact with at least one blinking neighbor do not advance, as if their phase update is being inhibited by the blinking neighbors, and otherwise advance to the next color (see Figure \ref{fig:FCA_ex1}). More precisely, the coupling for FCA is defined as the following equation  
\begin{align}\label{eq:def_FCA}
X_{t+1}(v)=\begin{cases} X_{t}(v)  & \text{if $X_{t}(v) > b(\kappa)$ and $X_{t}(u)=b(\kappa)$ for some $u\in N(v)$} \\
X_{t}(v)+1 & \text{otherwise.}\end{cases} 
\end{align}
For visualization purposes, it is convenient to consider the equivalent dynamics of `centered colorings' $\bar{X}_{t}:=X_{t}-t$ (mod $\kappa$), so that if $X_{t}$ synchronizes, then $\bar{X}_{t}$ converges to a constant function. In fact, FCA dynamics are displayed in this way in Figures 1 and 2 in the main text.

On the other hand, the \textit{Greenberg-Hastings model} (GHM) is a discrete model for neural networks \cite{greenberg1978spatial} introduced in 1987 by Greenberg and Hastings, where nodes of color $0$ and $1$ are called `rested' and `excited', respectively, and of any other color called 'refractory'. As in biological neural networks, rested neurons gets excited by neighboring excited nodes, and once excited, it has to spend some time in rested states to come back to the rested state again.  More precisely, the coupling for GHM is defined as   
\begin{align}\label{eq:def_GHM}
X_{t+1}(v) =
\begin{cases} 0  & \text{if $X_{t}(v)=0$ and $X_{t}(u)\ne 1$ for all $u\in N(v)$} \\ 
1  & \text{if $X_{t}(v)=0$ and $X_{t}(u)= 1$ for some $u\in N(v)$} \\
X_{t}(v) + 1 & \text{otherwise.}
\end{cases} 
\end{align}
For GHM, note that $X_{s}$ is synchronized if and only if $X_{t}\equiv 0$ for all $t\ge s+\kappa$.

In all experiments in this paper, we consider $\kappa=5$ instances of FCA and GHM. From here and hereafter, by FCA and GHM we mean the 5-color FCA and 5-color GHM, respectively. 

\begin{figure}
	\centering
	\vspace{-0.4cm}
	\includegraphics[width=0.6\textwidth]{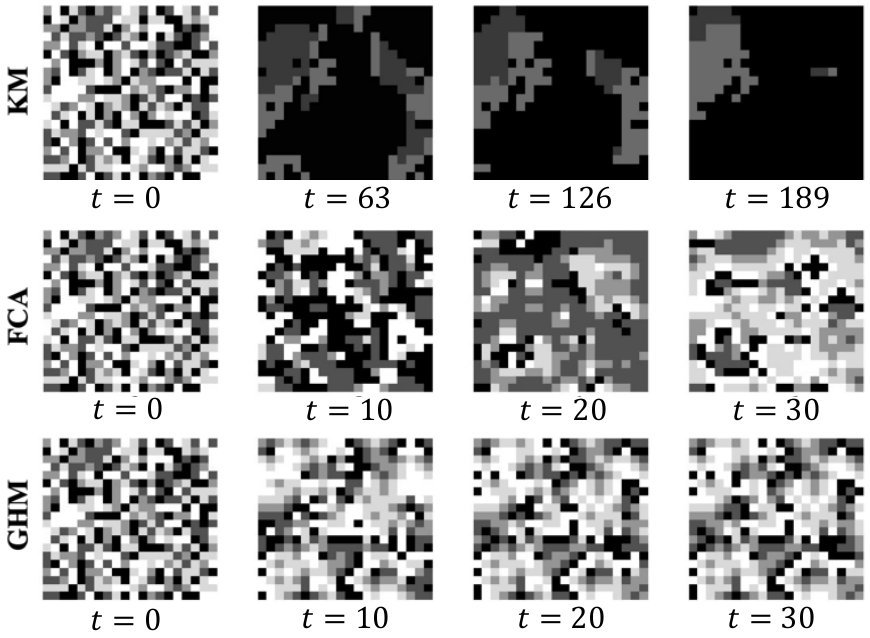}
	\caption{Simulation of KM, FCA, and GHM on the same underlying graph and initial configuration. The graph is generated by adding 80 edges uniformly at random into $20\times 20$ square lattice.  Each square heap map represents a phase configuration $X_{t}$ at the corresponding iteration $t$ shown below.  
	} 
	\label{fig:sim}
\end{figure}

\begin{figure}[H]
	\centering
	\vspace{-0.4cm}
	\includegraphics[width=1\textwidth]{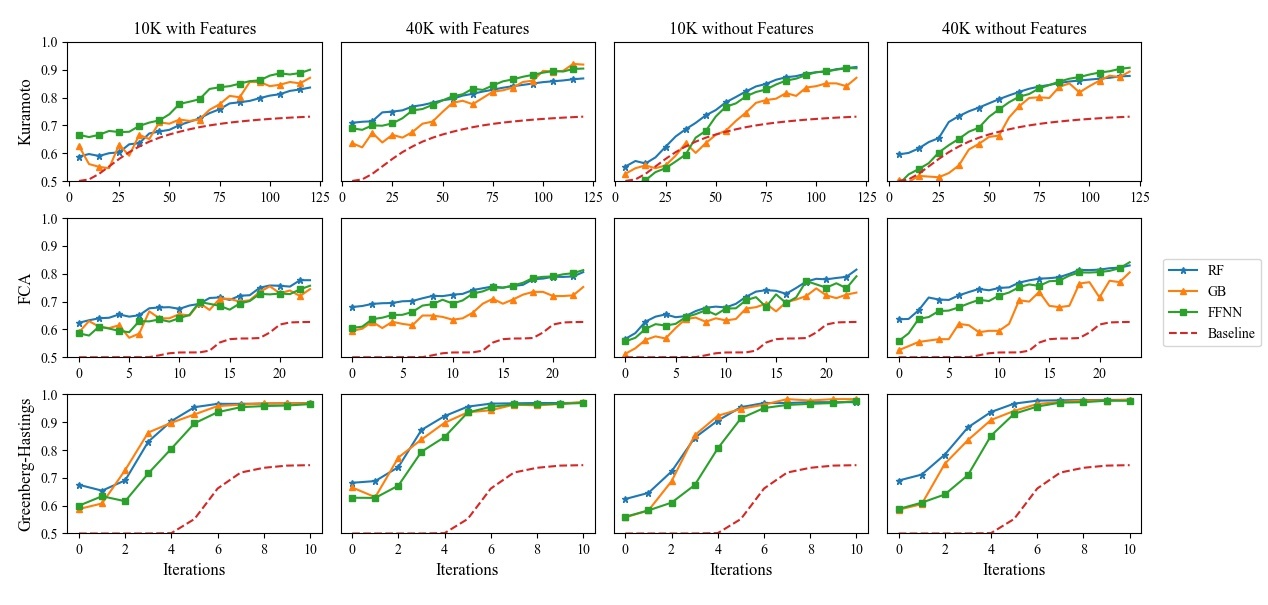}
	\caption{Supplementary experiment to Figure 5 in the main text to show sensitivity of our synchronization prediction methods to the size of training dataset. All settings are identical to the one in Figure 3 in the main text, except here we use smaller subsampled training set of sizes 10K and 40K for each classes of synchronizing and non-synchronizing examples, averaged over five runs. \commHL{The first two columns use datasets with 15 node graphs}, while the last two columns use datasets with 30 node graphs. Note that the full dataset used in Figure 3 in the main text uses 80K examples. Even trained with these smaller training set, we observe that most of our prediction algorithms significantly outperforms the baseline, although the performance does decrease in comparison to the main experiment in Figure 5 in the main text. (Note: This was averaged across 5 runs to address the randomness associated with subsampling). 
	} 
	\label{fig:volume}
\end{figure}

\begin{figure*}
    \centering
     \includegraphics[width=1\linewidth]{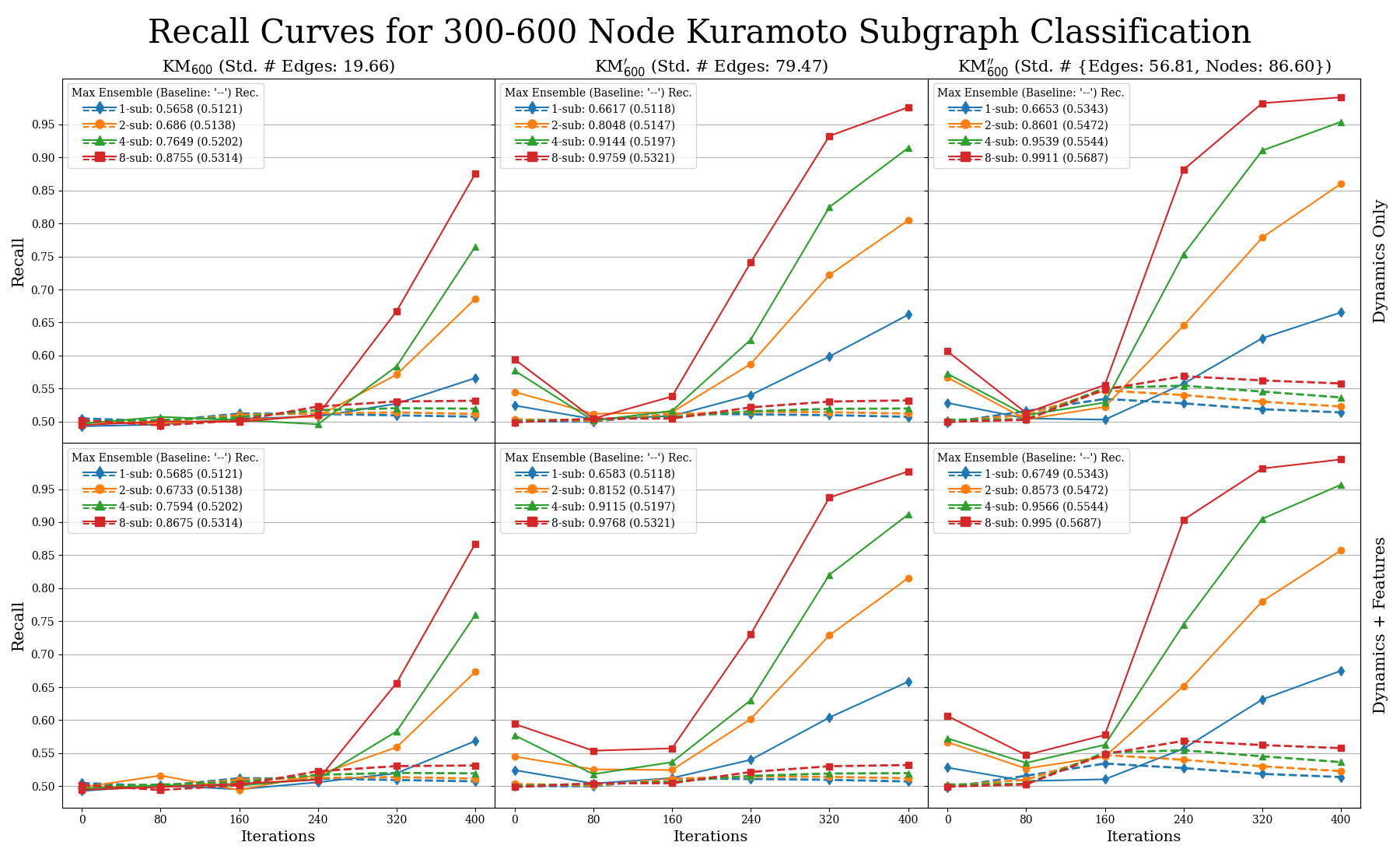}
    \centering
     \includegraphics[width=1\linewidth]{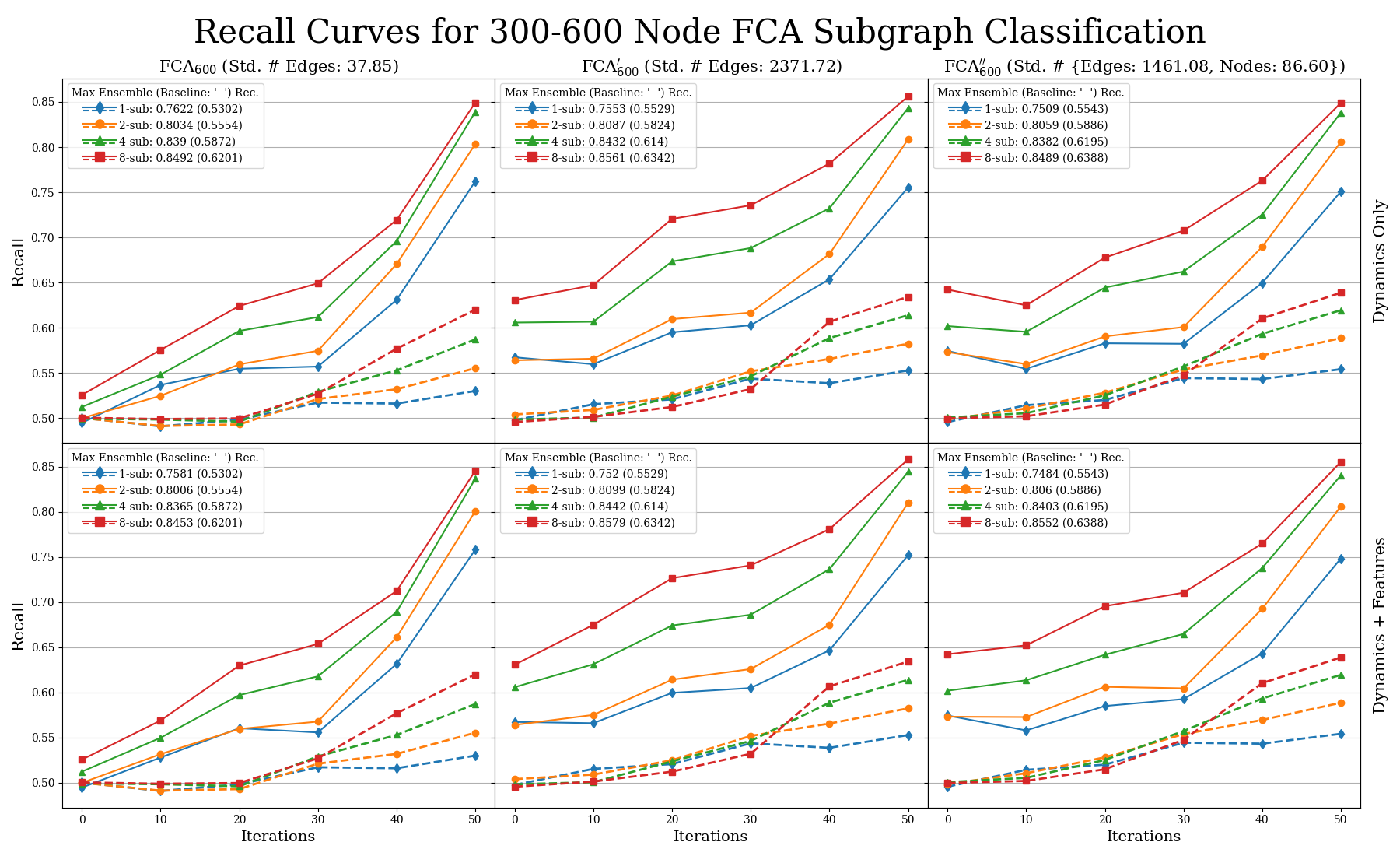}
\vspace{-0.3cm}
\caption{
 Recall curves for predicting synchronization of both Kuramoto (top) and 5-color FCA (bottom) on 600-node graphs from dynamics observed from $k\in\{1,2,4,8\}$ subgraphs of 30 nodes. All plots observe the performance of both the ensemble machine learning (solid) and baseline (dashed) recall scores over increasing amounts of training iterations  $r\in\{0,80,240,320,400\}$. The first row shows results using only dynamics whereas the second row includes both the dynamics and graph features. Maximum recall scores for using $k$ subgraphs are given by `$k$-sub: Acc. (Baseline Acc.)'
}
\label{fig:ensemble_plot_accuracy3}
\end{figure*}

\begin{figure*}[htpb]
    \centering
     \includegraphics[width=1\linewidth]{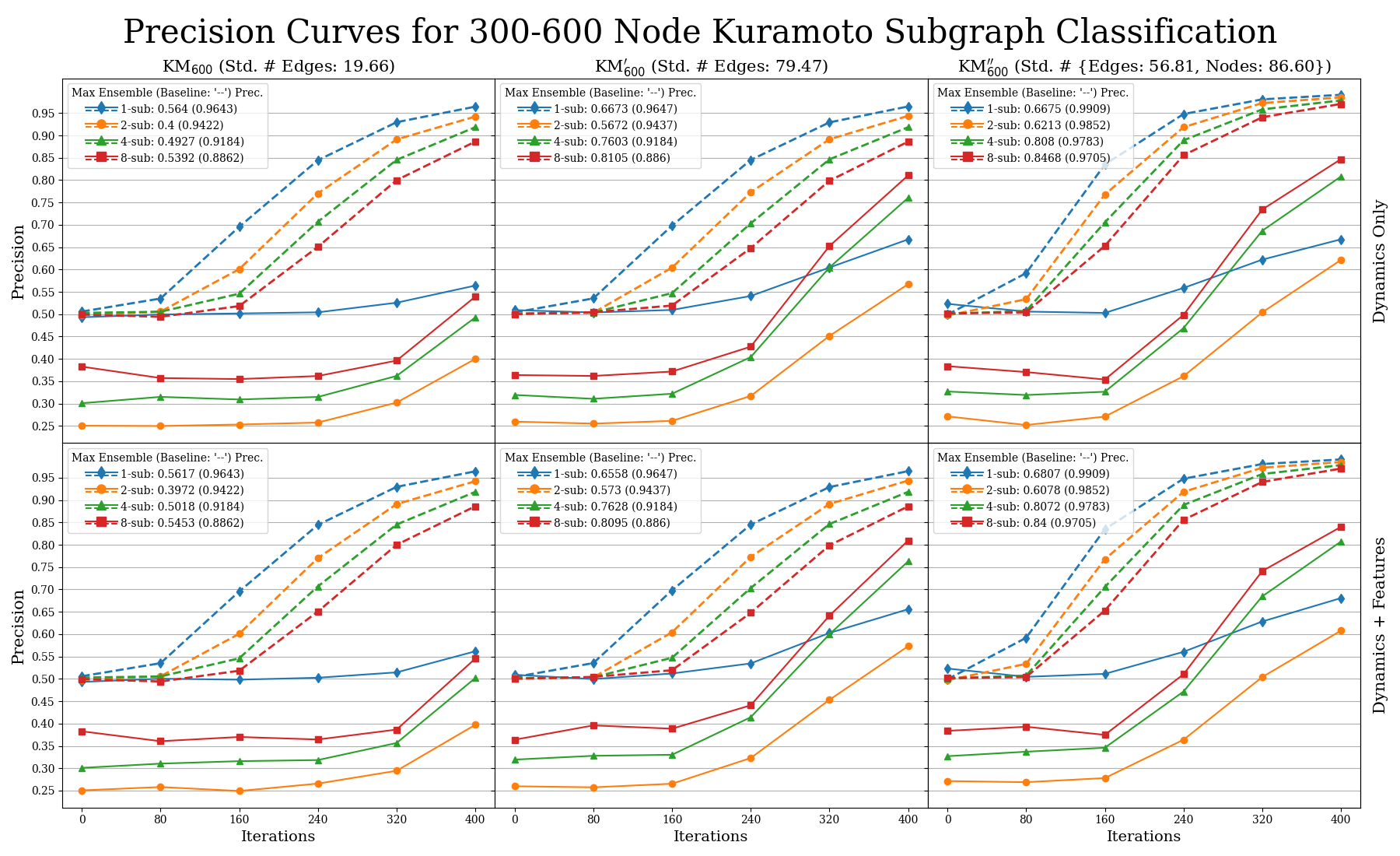}
    \centering
     \includegraphics[width=1\linewidth]{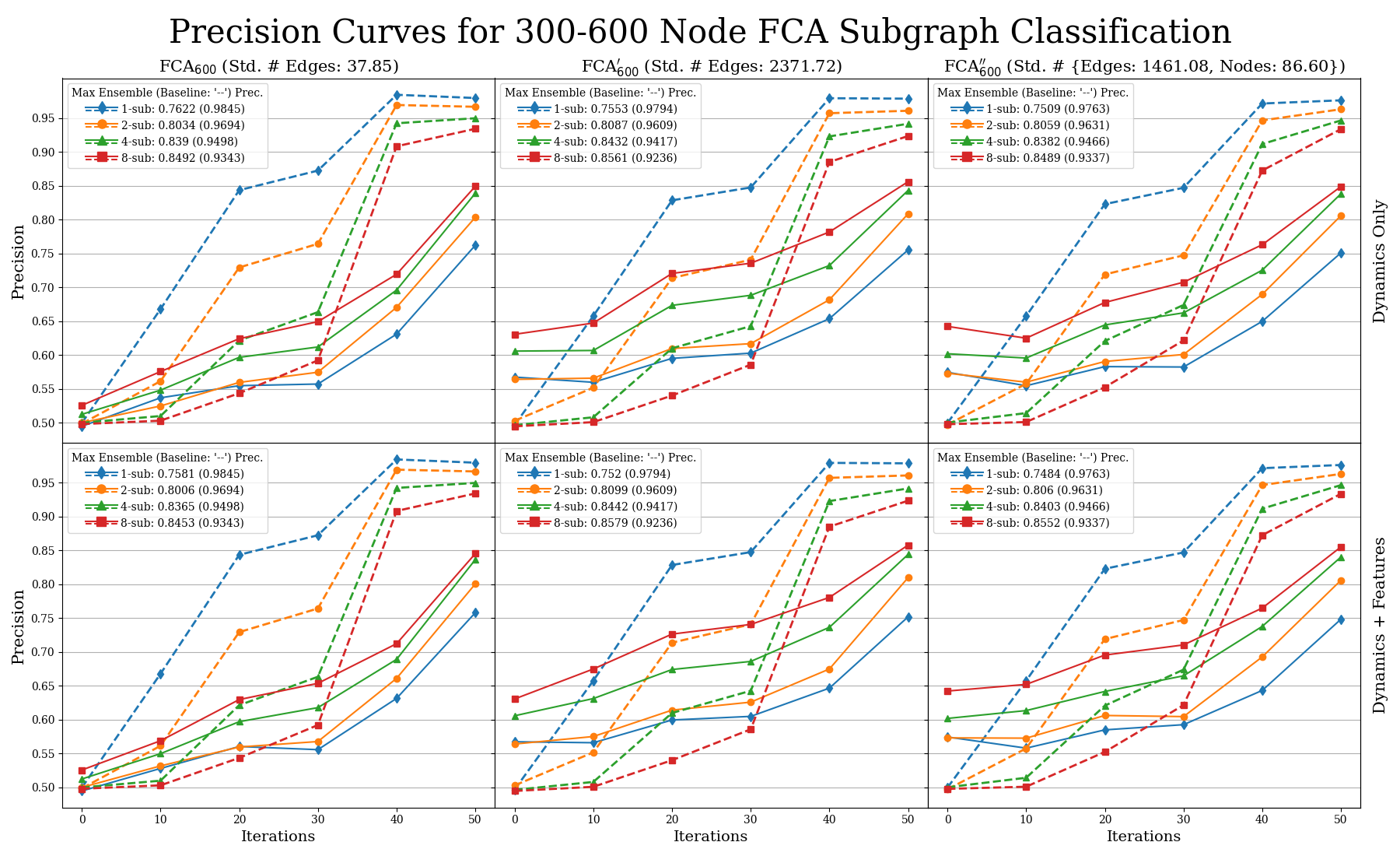}
\vspace{-0.3cm}
\caption{ Precision curves for predicting synchronization of both Kuramoto (top) and 5-color FCA (bottom) on 600-node graphs from dynamics observed from $k\in\{1,2,4,8\}$ subgraphs of 30 nodes. All plots observe the performance of both the ensemble machine learning (solid) and baseline (dashed) precision scores over increasing amounts of training iterations  $r\in\{0,80,240,320,400\}$. The first row shows results using only dynamics whereas the second row includes both the dynamics and graph features. Maximum precision scores for using $k$ subgraphs are given by `$k$-sub: Acc. (Baseline Acc.)'}
\label{fig:ensemble_plot_accuracy4}
\end{figure*}

While all three models tend to synchronize locally their global behavior on the same graph and initial configuration evolve quite differently, as seen in Figure \ref{fig:sim}. There, we simulate each system on the same graph, a 20x20 lattice with an additional 80 edges added uniformly at random. A single initial phase concentration $X_{0}$ is chosen by assigning each node with a uniformly randomly chosen state from $\{0,1,2,3,4\}$. This initial configuration is evolved through three different models: FCA, GHM and KM. 

This suggests that it is not feasible to predict synchronization for all three dynamics in the same way. Furthermore, the inclusion of random edges to the square lattice may disrupt some well-known behavior (e.g., 
spiral waves for coupled oscillators in 2D lattice \cite{martens2010solvable}) and result in unpredictable dynamics, analytically so. Hence, predicting synchronization on fully heterogeneous sets of graphs is a difficult task since keeping track of local interactions of oscillators is challenging to do when the overall structures of graphs within the heterogeneous set correspondingly have highly irregular and diverse couplings.

\subsection*{A simple lemma about Greenberg-Hastings model on a ring}

\vspace{0.2cm}
\textbf{Lemma.} \textit{Consider the $\kappa$-color  Greenberg-Hastings model \eqref{eq:def_GHM} on a path of $n$ nodes. Then for arbitrary initial coloring, $\kappa\ge 3$, and $n\ge 1$, the system will always converge to an all-zero configuration in $n+\kappa$ iterations. }

\vspace{0.2cm}

\noindent \textit{Proof.} We use a back-tracking argument that was used to analyze FCA in \cite{lyu2015synchronization, lyu2016phase, lyu2019persistence}. Order the nodes in the path as $1,\dots,n$ from left to right so that nodes that differ by 1 are adjacent.  Let $X_{t}:\{1\dots,n\}\rightarrow \{0,1,\dots,\kappa\}$ denote the color configuration at time $t$. Fix $T>\kappa$  and suppose that $X_{T}$ is not identically zero. Then necessarily $X_{T-\kappa}$ should have at least one 1. 

Suppose that node $i$ has color 1 at time $T-\kappa$. Then by the dynamics, one of its two neighbors ($i-1$ or $i+1$) should have color 1 at time $T-\kappa-1$. Suppose that node $i+1$ has color $1$ at time $T-\kappa-1$. By a similar argument, one of its two neighbors, now $i$ or $i+2$, should have color 1 at time $T-\kappa-2$. But since $X_{T-\kappa-1}(i)=1$, we should have $X_{T-\kappa-2}(i)=0$, so we must have $X_{T-\kappa-1}(i+2)=1$. Repeating this argument, we see that backtracking in time, the color 1 at node 0 at time $T-\kappa$ should travel at constant unit velocity, that is, $X_{T-\kappa-s}(i+1+s)=1$ for all $s$ such that $i+1+s\le n$. It follows that $X_{T-\kappa-(n-i-1)}(n-1)=0$ and $X_{T-\kappa-(n-i-1)}(n)=1$. But this yields a contradiction since then at time $T-\kappa-(n-i)$, node $n$ has color 0 and its only neighbor $n-1$ cannot have color 1 to excite node $n$. By a similar argument, if node $i-1$ has color $1$ at time $T-\kappa-1$, then the dynamics can only be back-tracked through time $T-\kappa-i+1$. This shows that $X_{T-\kappa}$ cannot have color 1 if $T-\kappa>n$. We conclude that that $X_{T}$ should be identically zero of $T>\kappa + n$. \hfill $\square$.

\subsection*{Machine learning algorithms for binary classification}

Different machine learning algorithms were employed to solve the problem of predicting whether or not a given graph and initial coloring for a system of coupled oscillators will tend to a synchronizing or non-synchronizing trajectory.

\begin{description}
    \item[Random Forest (RF) \cite{breiman2001random} ] A random forest (RF) is a form of ensemble learning that produces decision trees, and in which these trees themselves are generated using random feature vectors that are sampled independently from a distribution.
    
    Our implementation imposes a limit on the maximum amount of features used per tree to be the square root of the total amount of features in our dataset, $\sqrt p$ (where $p$ is the number of features), and uses 100 estimator trees with iterates terminating at complete node purity.
    
    \item[Gradient Boosting (GB) \cite{friedman2002stochastic}] 
  
    Gradient boosting (GB) is an ensemble learning algorithm for classification tasks similar to RFs where a loss function is minimized by a collection of decision trees to form a strong learner from a group of weak learners. The main difference is that the trees are not trained simultaneously, as in RFs, but added iteratively, so that the current loss is reduced cumulatively so in each iteration for each tree. 
    
    Our implementation uses a learning rate of 0.4, 100 estimator trees, and the square root of total number of features per tree when searching for splits.
    
    \item[Feed-forward Neural Networks (FFNN) \cite{bishop2006pattern}] 
    
Feedforward Neural Networks (FFNN) form a class of function approximators consisting of multiple layers of connected linear and non-linear activation functions \cite{csaji2001approximation, hornik1991approximation}. The linear maps are parameterized by `weights', which are subject to training on a given dataset. The algorithm works by iterating two phases. The first phase is the `feed-forward' phase; the input is mapped to an output by applying the layers of linear and non-linear activation functions with the current weights. The loss between the output and the target labels are then computed. The second phase is the `back-propagation' phase; the weights are modified in the direction that reduces the computed loss from the forward propagation phase.

Our implementation of an FFNN uses a learning rate of 0.01, cross-entropy loss, 4 fully connected linear layers, a hidden layer size of 100, ReLU activation, batch normalization, batch size 256, and dropout of 0.25 probability across 35 epochs. See \cite{svozil1997introduction} for backgrounds and jargons. 
    
\item[GraphLRCN]  

A Long-term Recurrent Convolutional Network (LRCN) \cite{donahue2015long} is a neural network architecture that has been developed for video data classification by combining convolutional neural networks \cite{krizhevsky2012imagenet} and a special type of recurrent neural network known as the Long Short-Term Memory network (LSTM-Net)\cite{hochreiter1997long}. We propose a variant of LRCN that is suitable for learning to predict dynamics on graphs, which we call the \textit{GraphLRCN}. The idea is to encode each configuration $X_{t}$ on a graph $G$ as a weighted adjacency matrix $\Delta(X_{t})$ of $G$. Then the dynamics $(X_{t})_{0\le t\le r}$ can be turned into a sequence of square matrices, which can be viewed as video data subject to classification.

Here we give a precise definition of the encoding $X_{t}\mapsto \Delta(X_{t})$. Given a $\kappa$-coloring $X$ on a graph $G=(V,E)$, we define the associated \textit{color displacement matrix} $\Delta=\Delta(X)\in (\mathbb{Z}_{\kappa})^{|V|\times |V|}$ as 
\begin{align}\label{eq:Delta_encoding}
    \Delta(i,j) = \min\left( \begin{matrix} X(i)- X(j) \,\, (\textup{mod $\kappa$}), \\ X(j)-X(i) (\textup{mod $\kappa$}) \end{matrix}  \right) \mathbf{1}(\text{$(i,j)\in E$}).
\end{align}
If $X:V\rightarrow [0,2\pi)$ is a configuration for KM, we define the associated color displacement matrix $\Delta(X)$ similarly as above by replacing $\textup{mod $\kappa$}$ with \textup{mod $2\pi$}. One can think of $\Delta(X)$ as the adjacency matrix of $G$ weighted by the color differences assigned by $X$ on the edges. An important feature of this encoding is that $\Delta(X)(i,j)$ is nonzero if and only if  nodes $i,j$ are adjacent in $G$ \textit{and} have different colors $X(i)\ne X(j)$. For instance, if $X$ is synchronized then $\Delta(X)$ is the zero matrix. Additionally, since the convolutional block components learn physical location-based associations in the matrix, we applied the reverse Cuthill-Mckee algorithm\cite{cuthill1969reducing} on the sequences of adjacency matrices to reduce our matrix bandwidths and augment these physical associations in our representations. In combination, the techniques mentioned are intended to reduce the amount of unnecessary information and noise in the encoding for the purpose of learning to predict synchronization. 
Our implementation uses a learning rate of 0.01, a batch size of $2^{11}$, and batch normalization across 25 epochs.
\end{description}

\subsection*{Additional figures}

\label{subsection:add_fig}

\end{document}